\documentstyle[11pt,epsfig]{article}

\newtheorem{The}{Theorem}[section]
\newtheorem{Cor}[The]{Corollary}
\newtheorem{Pro}[The]{Proposition}

\newtheorem{Rem}[The]{Remark}
\newtheorem{Lem}[The]{Lemma}

\def\proof{\vspace{2ex}\noindent{\bf Proof.} }

\def\endproof{\noindent $\diamond$}

\def\dirac{\partial}
\def\d{\partial}

\newcommand{\C}{{\bf C}}
\newcommand{\R}{{\bf R}}
\newcommand{\Z}{{\bf Z}}

\newcommand{\A}{{\cal A}}
\newcommand{\E}{{\cal E}}
\newcommand{\M}{{\cal M}}
\newcommand{\D}{{\cal D}}

\newcommand{\la}{\langle}
\newcommand{\ra}{\rangle}
\newcommand{\ba}{\begin{eqnarray}}
\newcommand{\na}{\end{eqnarray}}

\newcommand{\s}{{\bf s}}
\newcommand{\spinc}{\mathrm{Spin}^c}

\title{Exact triangles in Seiberg-Witten Floer theory. Part II: geometric
limits of flow lines} 
\author{Matilde Marcolli, Bai-Ling Wang}
\date{}

\begin{document}
\maketitle

\tableofcontents

\section{Introduction}

This work is the continuation of \cite{CMW}. In the previous part we
analyzed the geometric limits of solutions of the Seiberg--Witten
equations on a 3-manifold $Y(r)$ with a long cylinder $T^2 \times
[-r,r]$. We applied the result to the case of a homology sphere $Y$
with an embedded knot $K$ along which Dehn surgery with framing one or
zero is performed to
produce a homology sphere $Y_1$ or a manifold $Y_0$ with the homology of
$S^1\times S^2$. We proved that there is a suitable choice of
metrics on $Y$ and $Y_1$, and of a ``surgery perturbation'' $\mu$ of
the equations on a tubular  
neighborhood of the knot in $Y$, such that the moduli spaces of
solutions on $Y$, $Y_1$, and $Y_0$ are  related by 
$$ {\cal M}_{Y,\mu} \cong {\cal M}_{Y_1} \cup \bigcup_{k} {\cal
M}_{Y_0}(\s_{k}), $$
where the $\s_{k}$ are the possible choices of $\spinc$ structures
on $Y_0$. We also analyzed the splitting of the spectral flow and the
relative grading of the Floer complexes on $Y$, $Y_1$, and $Y_0$.

In this paper we concentrate on the analysis of
the geometric limits of flow lines, that is, of solutions of the
Seiberg--Witten equations on the four-manifold $Y(r)\times \R$, as
$r\to\infty$. This is necessary for an understanding of how the
boundary operators in the Floer complexes of $Y$, $Y_1$, and $Y_0$ are 
related. 

We first consider the unperturbed equations on $Y(r)\times \R$. We
derive estimates that control the convergence on compact sets. The
essential phenomenon which regulates the behavior of solutions is a
non-uniformity of the  
convergence in the non-compact time direction $t\in\R$, as we stretch
the length of the cylinder $r\to\infty$. This 
can be seen as an effect of the presence of small eigenvalues of the
linearization of the Seiberg-Witten equations on $Y(r)$ at the two
asymptotic values corresponding to $t\to \pm\infty$. In fact, as
already observed in 
\cite{CMW}, the linearization at these asymptotic values is a
self-adjoint first order elliptic operator with small
eigenvalues decaying like $1/r$, as $r\to\infty$. 

We can decompose the manifold $Y$ as $Y=V\cup_{T^2} \nu(K)$, where $V$
is the knot complement and $\nu(K)$ is a tubular neighborhood of the
knot. We use the same notation $V$ and $\nu(K)$ in the following to
indicate the manifolds completed with an infinite cylindrical end of
the form $T^2\times [0,\infty)$. 

We first analyze the properties of the geometric limits obtained in
the convergence of solutions $(\A_r,\Psi_r)$ on fixed compact sets
independent of $r\geq r_0$.
We show that on the two sides $V\times \R$ and $\nu(K)\times \R$ these
limits on compact sets decay exponentially in radial gauge  
along the region $T^2\times [0,\infty)\times \R$ to an asymptotic value
that is (up to gauge) a constant flat connection on $T^2$. 

We then analyze the convergence in asymptotic regions of the form
$Y(r)\times ([T_r,\infty)\cup(-\infty,-T_r])$, ``away from compact
sets''. To this purpose, we introduce a suitable rescaling of the
coordinates in such a way as to report within a finite region the
behavior in the complement of an arbitrarily large compact
sets. After the rescaling, we obtain convergence to paths in 
${\cal M}_V$ and ${\cal M}_{\nu(K)}$, and to holomorphic maps in
$\chi_0(T^2,V)$ and $\chi_0(T^2,\nu(K))$. 

These data match the geometric limits on the flat cylindrical region
$T^2\times [-r,r]\times \R$, which also consist of a flat connection,
as limit uniformly on compact set, and a holomorphic map as limit 
after the rescaling. Thus, we can
use the data of the geometric limits described above in order to form
approximate solutions on $Y(r)\times \R$, for $r(T,\tau)$ depending on 
the rescaling parameters, used here as gluing parameters. We prove
surjectivity of the linearization 
of the Seiberg--Witten equations at these approximate solutions. Thus, 
the approximate solutions can be deformed to actual solutions on
$Y(r)\times \R$. In other words, we obtain a complete analogue of the
gluing theorem for solutions of the 3-dimensional Seiberg--Witten
monopole equations analyzed in \cite{CMW}.

The rescaling technique and radial gauge limits that we consider in
this paper can be thought of as an analogue of the analysis of
non-abelian ASD equations and holomorphic curves of \cite{DoSa}.

{\bf Acknowledgments:} We are deeply grateful to Tom Mrowka who
corrected several mistakes in our understanding of the
convergence and gluing of flow lines, and referred us to the results
of \cite{DoSa} and to their rescaling technique as a possible model for our
problem. The first author 
benefited of several conversations with T.R. Ramadas, who
suggested the use of the radial gauge limits. Part of this work was
done during visits of the first author to the Tata Institute of
Fundamental Research in Mumbai, and to the Max Planck Institut
f\"ur Mathematik in Bonn. We thank these institutions for the kind
hospitality and for support. The first author is partially supported
by NSF grant DMS-9802480. The second author is supported by ARC
Fellowship. 

\section{Energy and curvature estimates}

Let $Y$ be a compact oriented smooth three-manifold with a fixed
trivialization of the tangent bundle, with a local basis $\{ e_i
\}$. This trivialization determines a $Spin$-structure on $Y$, with
spinor bundle $S$. Twisting $S$ with a line bundle $L$ gives a
$\spinc$-structure $\s$, with spinor bundle $W$ and determinant
$\det(W)=\det(S)\otimes L^2$. For the purpose of this paper we shall be
concerned with the case of a homology sphere $Y$, with the unique
choice of the trivial $\spinc$-structure, or with the case of $Y_0$, a
3-manifold 
with the homology of $S^1\times S^2$, with the infinite family of
$\spinc$ structures $\s_k$ with $c_1(\det(W_k))=2k$, as discussed in
\cite{CMW}. 

On the manifold $Y(r)$, consider a solution $(A_r,\psi_r)$ of the
Seiberg--Witten equations 
\ba \begin{array}{c} {}*F_A = \sigma(\psi,\psi) \\[2mm]
\dirac_A \psi =0, \end{array}\label{3SW} \na
given in terms of a $U(1)$-connection $A$ on the line bundle $L$, and
a section $\psi\in \Gamma(Y,W)$. The 1-form $\sigma(\psi,\psi)$ is
given in local coordinates by $\sigma(\psi,\psi)=\sum_i \la e_i .
\psi,\psi \ra e^i$.
Solutions are critical points of the Chern--Simons--Dirac functional
$$ CSD(A,\psi)=-\displaystyle{\frac {1}{2} \int_Y} (A-A_0) \wedge 
(F_A +F_{A_0} ) + \displaystyle{\int_Y} \la \psi, \dirac_A \psi\ra
dvol_Y. $$
For the moment, we only consider the unperturbed functional and
the unperturbed equations (\ref{3SW}). Later in this paper 
it will be necessary to formulate the results in the perturbed case,
with the perturbations introduced in \cite{CMW}.

Let $L_{A_r,\psi_r}$ be the linearization of the equations, namely the 
operator 
$$ L_{A_r,\psi_r}(\alpha,\phi)= \left( -*d\alpha +\sigma (\psi_r, \phi),
\dirac_{A_r} \phi + \frac 12 \alpha . \psi_r \right). $$
This, together with the infinitesimal action of the gauge group and
the gauge fixing condition, defines the extended Hessian of the
Chern-Simons-Dirac functional, namely the first order elliptic operator
$$ \begin{array}{cl}
H_{A_r,\psi_r}(f,\alpha,\phi)= & ( -*d\alpha +\sigma (\psi_r, \phi)
-df,\\[2mm] &  \dirac_{A_r} \phi + \frac 12 \alpha .\psi_r + 
f\psi_r,\\[2mm] & -d^*\alpha -iIm \la \psi_r, \phi\ra ). \end{array}$$

In \cite{CMW} we proved the following gluing theorem. We assigned a 
metric on $Y$ which has positive scalar curvature on a tubular
neighborhood $\nu(K)$ of the knot. We proved that, with this
choice of the metric, for large $r\geq R_0$, the 
solutions $(A_r,\psi_r)$ can be written as a gluing
$$ (A_r,\psi_r)=(A',\psi')\#_r (A'',0). $$
Here $(A',\psi')$ is a finite energy solution of (\ref{3SW}) on the
knot complement $V$ endowed with an infinite cylindrical end
$T^2\times [0,\infty)$, and $(A'',0)$ is a flat connection on
$\nu(K)$. The length of the cylinder $r\geq R_0$ appears as the gluing
parameter. The gluing happens at the asymptotic value $a_\infty$ of
$(A',\psi')$ along the end $T^2\times [0,\infty)$ of the knot
complement $V$: $a_\infty$ is a flat connection on $T^2$, which lives
in a cyclic cover $\chi_0(T^2,V)$ of the character variety
$\chi(T^2)$. Up to a gauge transformation, $a_\infty$ agrees with the
restriction to $T^2$ of the connection $A''$.

In this paper we formulate the analogous gluing theorem for the
gradient flow lines of the Chern--Simons--Dirac functional, that is,
for solutions of the Seiberg--Witten equations on $Y(r)\times \R$. The
gluing theorem, in this case, is more complicated because of the lack of
uniformity in the convergence as $r\to\infty$. This gives rise to
more complicated geometric limits.
We report here a statement from \cite{CMW} which is the
ultimate source of the interesting non-uniform limits of flow lines
that we are going to discuss in the rest of this paper. 

\begin{Pro}
The dimension $N(r, o(1))$ of the span of eigenvectors of
the operator $ H_{A_r,\psi_r}$ with eigenvalues
satisfying $\mu(r)\to 0$ as $r\to\infty$ is given by
$$ N(r, o(1))=\dim Ker_{L^2}(H_{(A',\psi')}) + \dim Ker_{L^2}(
H_{(A '',0)}) + \dim Ker(Q_{a_\infty}), $$
where $a_\infty$ is the asymptotic value of $(A',\psi')$ as
$s\to\infty$.

For any $\epsilon >0$, the dimension $N(r,r^{-(1+\epsilon)})$ of the
span of eigenvectors of the operator $ H_{A_r,\psi_r}$ with eigenvalues
$\mu < r^{-(1+\epsilon)}$ is given by
$$ N(r,r^{-(1+\epsilon)})= \dim Ker_{L^2}(H_{(A',\psi')}) + \dim Ker_{L^2}(
H_{(A '',0)}) + \dim \ell_1 \cap \ell_2, $$
where $\ell_1 \cap \ell_2$ is the intersection of the two Lagrangian
submanifolds determined by the
extended $L^2$ solutions of $H_{(A',\psi')}(\alpha,\phi)=0$ and
$H_{(A '',0)}(\alpha,\phi)=0$ in the tangent space
$H^1(T^2,i\R)=Ker(Q_{a_\infty})$ of $\chi_0(T^2)$.

Thus, under the assumption that 
$$ Ker_{L^2}(H_{(A',\psi')})=Ker_{L^2}(H_{(A '',0)}) =0, $$
we have $N(r,o(1))=2$, generated by the elements of $H^1(T^2,i\R)$,
and $N(r,r^{-(1+\epsilon)})=0$, for all $\epsilon >0$.
\label{CLM}
\end{Pro}

The proof of Proposition \ref{CLM} relies on the results of \cite{CLM}, cf.
\cite{CMW}. 

Let $Y(r)$ be a closed three-manifold with a long tube $[-r,r]\times T^2$,
endowed with the flat product metric on the tube.

Let $(\A_r,\Psi_r)$ be a finite energy solution of the
Seiberg-Witten equations on the four-manifold $Y(r)\times \R$, that
is, $(\A_r,\Psi_r)$ satisfies the equations
\ba
\begin{array}{c} F_{\A}^+ = \tau(\Psi,\Psi) \\[2mm]
D_{\A}\Psi=0, \end {array}\label{4SW}
\na
with the condition
$$ \E_r=\| \partial_t A_r \|^2_{L^2(Y(r)\times\R)}+ \| \partial_t \psi_r
\|^2_{L^2(Y(r)\times\R)} <\infty. $$
Here $(A_r(t),\psi_r(t))$ is a solution in a temporal gauge, gauge
equivalent to  
the original $(\A_r,\Psi_r)$. The self-dual form $\tau(\Psi,\Psi)$ is
given in local coordinates by $\tau(\Psi,\Psi)=\sum_{i,j} \la e_i e_j
\Psi,\Psi \ra e^i\wedge e^j$.
As proved in \cite{MW}, the finite energy condition ensures the 
existence of asymptotic values $(A_r(\pm\infty),\psi_r(\pm\infty))$ as 
$t\to\pm\infty$, for any fixed $r\geq r_0$. The asymptotic values satisfy 
the Seiberg-Witten equations (\ref{3SW}) on $Y(r)$.

If the asymptotic values are
irreducibles, satisfying  $\psi_r(\pm\infty)\neq 0$, then the convergence
as $t\to\pm\infty$ is exponential,
\ba
\label{exp:decay}
\| (A_r(t),\psi_r(t)) - (A_r(\pm\infty),\psi_r(\pm\infty))
\|_{L^2_1(Y(r)\times \{ t \} )} \leq C_r e^{-\delta_r |t|}, 
\na
for all $t\geq T_r$. The rate of decay $\delta_r$ is determined by the
absolute value of the first non-trivial eigenvalue of the linearization 
$L_{A_r(\pm\infty),\psi_r(\pm\infty)}$
of the monopole equations (\ref{3SW}) at
$(A_r(\pm\infty),\psi_r(\pm\infty))$, hence, by Proposition \ref{CLM},
the rate of decay to the asymptotic
values at $\pm\infty$ satisfies
$$ \delta_r = \frac{c}{r}. $$ 

We are interested in studying the boundary operator of the Floer
complex, hence we are only interested in flow lines that connect
irreducible critical points. In the case of equivariant Floer theory
\cite{MW} we need also consider flow lines connecting to the reducible
point, but in the framed configuration space the reducible solution is
also a smooth point and a non-degenerate critical point, hence we have
the same exponential decay of flow lines, as proved in \cite{MW}. 

\begin{Lem}
Let $(\A,\Psi)$ be a solution of the Seiberg-Witten equations
(\ref{4SW}) on
$Y(r)\times \R$, for a fixed $r>0$. Suppose that $(\A,\Psi)$ is in a
temporal gauge on $Y(r)\times [T,\infty)$ and $Y(r)\times (-\infty,
-T]$ and assume that it decays exponentially to 
asymptotic values $(A(\pm\infty),\psi(\pm\infty))$ as $t\to\pm\infty$,
where the $(A(\pm\infty),\psi(\pm\infty))$ satisfy the 
Seiberg-Witten equations (\ref{3SW}).

Then we have the following identity
$$ -\frac 12 \int_{Y(r)}(A(-\infty) -A_0) \wedge 
(F_{A(-\infty)}+F_{A_0})+ \frac 12 \int_{Y(r)}(A(\infty) -A_0) \wedge
(F_{A(\infty)}+F_{A_0}) $$
$$ ={\frac 12  \int_{Y(r)\times\R }} F_\A
\wedge F_\A. $$
The identity holds more generally when replacing $Y(r)$ with a smaller
domain $\Omega \subset Y(r)$.
\label{energy}
\end{Lem}

\proof
We have
$$ -\frac 12 \int_{Y(r)}(A(-\infty) -A_0) \wedge 
(F_{A(-\infty)}+F_{A_0})+ \frac 12 \int_{Y(r)}(A(\infty) -A_0) \wedge
(F_{A(\infty)}+F_{A_0}) $$
$$ \begin{array}{lll}
&= &- {\frac 12  \int_{Y(r)}} (A(-\infty) -A_0) \wedge 
(F_{A(-\infty)}+F_{A_0})\\[2mm]
&&    +  {\frac 12  \int_{Y(r)}}
 (A(\infty) -A_0) \wedge (F_{A(\infty)}+F_{A_0})\\[2mm]
&= & {\frac 12  \int_\R dt {\frac{\partial}{\partial t}} \int_{Y(r)} }
\bigl( (A(t)-A_0) \wedge (F_{A(t)} + F_{A_0})\bigr)\\[2mm]
&= &-  {\frac 12  \int_{Y(r)}\int_\R }
 {\frac {\partial}{\partial t}}
F_{A(t)} \wedge (A(t)-A_0)\wedge dt \\[2mm] 
&& - \frac 12  \int_{Y(r)}\int_\R
(F_{A(t)} + F_{A_0})\wedge {\frac {\partial}{\partial t}} A(t)
\wedge dt \\[2mm]
&= &-  {\frac 12  \int_{Y(r)\times \R} }
d_{Y(r)}(  {\frac {\partial}{\partial t}} A(t))
 \wedge (A(t)-A_0)\wedge  dt \\[2mm]
&& - \frac 12  \int_{Y(r)}\int_\R
(F_{A(t) }+ F_{A_0})\wedge {\frac {\partial}{\partial t}} A(t)
\wedge dt \\[2mm]
&= &-  {\frac 12  \int_{Y(r)\times \R}} 
(F_{A(t)}-F_{A_0})\wedge {\frac {\partial}{\partial t}} A(t)
\wedge dt\\[2mm]
&& - \frac 12  \int_{Y(r)\times \R}
(F_{A(t)} + F_{A_0})\wedge {\frac {\partial}{\partial t}} A(t)
\wedge dt \\[2mm]
&= &-   {\int_{Y(r)\times \R}}F_{A(t)}\wedge {\frac {\partial
A(t)}
{\partial t}}
\wedge dt\\[2mm]
&= &  {\frac 12  \int_{Y(r)\times\R }} F_\A \wedge F_\A.
\end{array}
$$
This completes the proof. 

\endproof

\begin{Lem}
Let $(\A_r,\Psi_r)$ be a solution of the Seiberg-Witten equations
(\ref{4SW}) on 
$Y(r)\times \R$, such that the gauge transformed element
$$\lambda_r(\A_r,\Psi_r)=(A_r(t),\psi_r(t))$$ 
in a temporal gauge
satisfies the finite energy condition
$$ \E_r=\| \partial_t A_r \|^2_{L^2(Y(r)\times\R)}+ \| \partial_t \psi_r
\|^2_{L^2(Y(r)\times\R)} = $$
$$ \int_{\R} \| \nabla CSD (A_r(t),\psi_r(t)) \|^2_{L^2(Y(r))} dt
<\infty. $$ 

Then, for any interval $[t_0,t_1]$ of length $\ell=t_1-t_0$, we have 
estimates 
$$ \| \Psi_r \|^4_{L^4(Y(r)\times [t_0,t_1])}\leq 8 \E_r + 2s_0^2
vol(Y(r_0))l, $$
$$ \int_{t_0}^{t_1} \| F_{A_r(t)} \|^2_{L^2(Y(r)\times \{ t \})} dt
\leq \E_r +s_0^2 vol(Y(r_0))l, $$
and
$$ \int_{t_0}^{t_1} \| \nabla_{A_r(t)}\psi_r(t) \|^2_{L^2(Y(r)\times
\{ t \})} dt \leq \E_r +s_0^2 vol(Y(r_0))l, $$
for all $r\geq r_0$.

Moreover, the energy $\E_r$ is uniformly bounded
in $r\geq r_0$. The constant $s_0\geq 0$ is defined by  
\ba\label{s0}
s_0 =\max_{Y(r_0)} \{ -s(x), 0 \}, \na 
with $s(x)$ the scalar curvature. 
\label{L1}
\end{Lem}

\proof
The argument follows 6.12 of \cite{MST}. Since the cylinder
$[-r,r]\times T^2$ is endowed with the flat metric, we have
$$ s(r)=\max_{Y(r)} \{ -s(x), 0 \} = s_0, $$
for $r\geq r_0$ and $s_0$ as in (\ref{s0}).
Moreover, we can estimate
$$ \begin{array}{c}
\E_r \geq \int_{t_0}^{t_1} (\frac{1}{2} \| F_{A_r(t)}
\|^2_{L^2(Y(r))} + \| \nabla_{A_r(t)} \psi_r(t) \|^2_{L^2(Y(r))}) dt + 
\\[2mm]
\int_{Y(r)\times [t_0,t_1]}\frac{1}{8} |\Psi_r|^4 dv_{Y(r)}dt
-\frac{s_0}{4}\int_{Y(r_0)\times [t_0,t_1]} |\Psi_r|^2 dv_{Y(r_0)}dt. 
\end{array} $$
At any local maximum, $|\Psi_r|^2$ is bounded by the scalar curvature,
hence a non-trivial maximum can only occur away from $[-r,r]\times
T^2$, i.e. on 
$Y(r_0)\times\R$. Thus, the estimate above gives
$$ \| \Psi_r \|^4_{L^4(Y(r)\times [t_0,t_1])}\leq 8 \E_r + 4s_0^2
vol(Y(r_0))l. $$
The other two estimates follow similarly.

The uniform bound on the energy $\E_r$ is obtained as follows. As we 
discussed before, for each fixed $r>0$, the finite energy condition for
$(\A_r,\Psi_r)$ forces the existence of asymptotic values
$(A_r(\pm\infty),\psi_r(\pm\infty))$ as $t\to\pm\infty$ that satisfy
the 3-dim Seiberg-Witten equations (\ref{3SW}). Moreover, if the
asymptotic values 
are non-degenerate critical points of the Chern-Simons-Dirac
functional, then the temporal gauge representative $(A_r(t),\psi_r(t))$
decays sufficiently fast in $t$ to the asymptotic values in the $L^2_2$
topology. 

Thus, the energy $\E_r$ can be written also as the total variation of the
Chern-Simons-Dirac functional along the path $(A_r(t),\psi_r(t))$,
$$ \E_r= CSD(A_r(-\infty),\psi_r(-\infty))-
CSD(A_r(\infty),\psi_r(\infty)). $$ 
Moreover, since the elements $(A_r(\pm\infty),\psi_r(\pm\infty))$
satisfy the equations (\ref{3SW}) on $Y(r)$, we have
$\dirac_{A_r(\pm\infty)} \psi_r(\pm\infty)=0$ and the variation of the CSD
functional is simply given by
$$ \E_r=\frac{1}{2}\int_{Y(r)} (A_r(\infty)-A_0)\wedge
(F_{A_r(\infty)}+F_{A_0}) $$
$$-\frac{1}{2} \int_{Y(r)} (A_r(-\infty)-A_0)\wedge
(F_{A_r(-\infty)}+F_{A_0})$$

The previous Lemma shows that this quantity can be rewritten as
$$ \E_r=\frac{1}{2} \int_{Y(r)\times  \R} F_{\A_r}\wedge F_{\A_r}. $$ 
This is a topological term, conformal and gauge invariant, hence it
does not change when stretching the cylinder $[-r,r]\times T^2$.

\endproof

\begin{Lem}
Let $(A_r(t),\psi_r(t))$ be a finite energy solution of the 
Seiberg-Witten equations (\ref{4SW}) on $Y(r)\times \R$, in temporal
gauge in the $\R$-direction.
Then, for all $t\geq T_r$, we can write the estimate (\ref{exp:decay}) 
as
$$ \| (A_r(t),\psi_r(t)) - (A_r(\pm\infty),\psi_r(\pm\infty))
\|_{L^2(Y(r)\times \{ t \} )} \leq C e^{-\frac{c|t|}{r}}, $$
where the constant $C$ is independent of $r$.  
Moreover, on any fixed compact set $K \subset Y(r)$, independent of $r\geq
r_0$, we have
$$ \| (A_r(t),\psi_r(t)) - (A_r(\pm\infty),\psi_r(\pm\infty))
\|_{L^2_1( K' \times \{ t \} )} \leq C_K e^{-\frac{c|t|}{r}}, $$
with $K' \subset int(K)$, and with $C_K$ independent of $r$.
\label{decay:bound}
\end{Lem}

\proof
The $L^2$ bound follows by the exponential decay of the energy
functional: for $t> T_r$ we have 
$$ \frac{d}{dt}\E_r(t) \leq \E_r(0) \exp(-\frac{ct}{r}), $$
where
$$ \E_r(t)=\int_t^\infty \| \partial_t A_r \|_{L^2(Y(r))}^2 + \|
\partial_t \psi_r \|^2_{L^2(Y(r))} dt, $$
and $\E_r(0)$ uniformly bounded in $r\geq r_0$, by Lemma \ref{L1}. 

The $L^2_1$ bound follows by retracing the argument of Lemma 6.14 of
\cite{MST}: we obtain that the constant $C_r$ of our (\ref{exp:decay})
depends on the underlying manifold $Y(r)$ through the constant of the
Sobolev multiplication theorem, and the estimates of our Lemma
\ref{L1}, which are uniform in $r\geq r_0$. Consider a compact set $K$
that is embedded in $Y(r)$ for all $r\geq r_0$, so that the metric on
$K$ does not change with $r\geq r_0$. If we restrict our estimate
to $K$, and we choose the cutoff function $\xi$ of
Lemma 6.14 of \cite{MST} supported in $K$, we obtain a uniform $C_K$
depending only on $K$.

\endproof

\begin{Lem}
Suppose given $(\A_r,\Psi_r)$, a finite energy solution of
the  equations (\ref{4SW}) on $Y(r)\times \R$, for $r\geq r_0$. 
There is a pointwise bound
$$ | \Psi_r (x,t) |^2 \leq s_0, $$
where $s_0 =\max_{Y(r_0)} \{ -s(x), 0 \}$, with $s(x)$ the scalar
curvature. Moreover, for any $t_0 < t_1$, we have an estimate 
$$ \int_{Y(r)\times [t_0,t_1]} \la \Psi_r,
\nabla^*_{\A_r}\nabla_{\A_r}\Psi_r \ra dv_{Y(r)}dt \leq
\frac{s_0}{4} \| \Psi_r \|^2_{L^2(Y(r_0) \times [t_0,t_1])}. $$
\label{Lpoint}
\end{Lem}

\proof
The finite energy condition ensures the existence of asymptotic values 
$\psi_r( x,\pm\infty)$ for the spinor $\Psi_r (x,t)$. Thus, we can
estimate 
$$ | \Psi_r (x,t) |^2 \leq \max \{ |\psi_r( x,\pm\infty) |^2,
-s_r(x,t) \}, $$
with $s_r(x,t)$ the scalar curvature on $Y(r)\times \R$,
and the pointwise estimate follows.
By the Weitzenb\"ock formula and the equations we have
$$ 0 = \nabla^*_{\A_r}\nabla_{\A_r} \Psi_r + \frac{s_r}{4} \Psi_r +
\frac{| \Psi_r |^2}{4} \Psi_r. $$
By integrating and by using the vanishing of the scalar curvature on
the cylinder $T^2\times [-r,r]$, we obtain the desired estimate.

\endproof

In the following we use the notation
$$ \E_r (t_0,t_1)=\int_{t_0}^{t_1} \| \nabla CSD (A_r(t),\psi_r(t))
\|^2_{L^2(Y(r))} dt. $$

\begin{Pro}
Assume, as before, that $(\A_r,\Psi_r)$ is a finite energy solution of
the equations (\ref{4SW}) on $Y(r)\times \R$. Suppose given any two real
numbers $t_0<t_1$. Then, there is a uniform bound on  
$\| F_{\A_r} \|_{L^2(Y(r)\times [t_0,t_1])}$ for $r\geq r_0$,
$$ \| F_{\A_r} \|^2_{L^2(Y(r)\times [t_0,t_1])} \leq C_\ell, $$
where $C_\ell$ is a constant independent 
of $r$ and depending only on $\ell=t_1-t_0$.
\end{Pro}

\proof
According to the previous Lemma we have a uniform bound on the $L^4$-norm
$$ \| \Psi_r
\|_{L^4(Y(r)\times [t_0,t_1])}, $$ 
hence, by the equations, we obtain a
uniform bound on $\| F_{\A_r}^+ \|^2_{L^2(Y(r)\times [t_0,t_1])}$. 

We have the identity 
\begin{equation}
\label{curv-energy}
 (2 | F_{\A}^+ |^2-| F_{\A} |^2) dv = - F_{\A}\wedge
F_{\A} .
\end{equation}
Thus, we obtain an estimate
$$ \begin{array}{ll}
\| F_{\A_r} \|^2_{L^2(Y(r)\times [t_0,t_1])} \leq 2  & \| F_{\A}^+
\|^2_{L^2(Y(r)\times [t_0,t_1])} + \E_r (t_0,t_1) \\[2mm]
 & + \int_{Y(r)}\la \psi_r(t_1),\partial_{A_r(t_1)}\psi_r(t_1) \ra
dv_{Y(r)} \\[2mm]
& - \int_{Y(r)}\la \psi_r(t_0),\partial_{A_r(t_0)}\psi_r(t_0) \ra
dv_{Y(r)}.
\end{array} $$
In the right hand side we can estimate 
$$ \E_r (t_0,t_1) \leq \E_r, $$
which is uniformly bounded in $r\geq r_0$, by Lemma \ref{L1}. We need
to estimate the term
$$ | \int_{Y(r)}\la \psi_r(t_1),\partial_{A_r(t_1)}\psi_r(t_1) \ra
dv_{Y(r)}  - \int_{Y(r)}\la \psi_r(t_0),\partial_{A_r(t_0)}\psi_r(t_0) \ra
dv_{Y(r)} |. $$

We estimate
$$ \begin{array}{c}
| \int_{t_0}^{t_1} \frac{d}{dt} \int_{Y(r)} \la
\psi_r(t),\partial_{A_r(t)}\psi_r(t) \ra dv_{Y(r)} dt | \leq \\[2mm] 
\int_{t_0}^{t_1} ( 2\| \partial_{A_r(t)}\psi_r(t) \|^2_{L^2(Y(r))} +
2 | \la F_{A_r(t)} , *\sigma(\psi_r(t),\psi_r(t)) \ra |) dt \leq \\[2mm]
\int_{t_0}^{t_1}  2\| \nabla_{A_r(t)}\psi_r(t) \|^2_{L^2(Y(r))} + \\[2mm]
( \int_{t_0}^{t_1} \| F_{A_r(t)} \|^2_{L^2(Y(r))} dt )^{1/2}  
( \int_{t_0}^{t_1} \int_{Y(r)} \| \psi_r(t) \|^4 dv_{Y(r)} dt )^{1/2} +
\\[2mm]
\frac{s_0}{2} \| \Psi_r \|^2_{L^2(Y(r_0)\times [t_0,t_1])} +
\frac{1}{2} \| \Psi_r \|^4_{L^4(Y(r)\times [t_0,t_1])} .\end{array} $$

Now the estimates of Lemma \ref{L1} give the required bound on $\| F_{\A_r}
\|$. 

\endproof

\section{Convergence on compact sets}

In the process of stretching the neck in $Y(r)\times \R$,
the spinor tends to vanish pointwise on the long cylinder, as the following 
Lemma proves.

\begin{Lem}
Assume, as before, that $(\A_r,\Psi_r)$ is a finite energy solution of
the equations (\ref{4SW}) on $Y(r)\times \R$. Given any compact set of
the form $T^2\times [-r_0,r_0]\times [t_0,t_1]$,  
we have pointwise convergence
$$ | \Psi_r(x,t) | \to 0, $$
for all $(x,t)\in T^2\times [-r_0,r_0]\times [t_0,t_1]$, as
$r\to\infty$. 
\label{pointwise:conv}
\end{Lem}

\proof
The claim follows from the uniform bound
$$ \int_{T^2\times [-r,r]\times [t_0,t_1]} | \Psi_r(x,t) |^4 dv dt
\leq C_\ell $$
derived previously.

\endproof

Thus, the self-dual part of the curvature is also converging to zero
pointwise on the long cylinder.

\begin{Cor}
If $(\A_r,\Psi_r)$ is a finite energy solution of
the equations (\ref{4SW}) on $Y(r)\times \R$, then, on compact sets of
the form $T^2\times [-r_0,r_0]\times [t_0,t_1]$, the connection $\A_r$
converges to a finite energy solution $\A$ of the abelian self dual
equation $F_{\A}^+=0$. 
\end{Cor}

\begin{Pro}
Inside the manifold $Y(r)\times \R$, for all $r\geq r_0$, consider  
fixed compact sets of the form $K_V \subset V\times \R$, $K_\nu
\subset \nu(K)\times \R$, or $K_0 \subset T^2\times [-r_0,r_0]\times
\R$. We assume that the metric on $\nu(K)$ has non-negative scalar curvature,
and that it is flat on the cylinder $T^2\times [-r,r]$. 
Let $(\A_r,\Psi_r)$ be finite energy solutions of the 
equations (\ref{4SW}) on $Y(r)\times \R$. Then there exists a subsequence
$(\A_{r'},\Psi_{r'})$ and gauge transformations $\lambda_{r'}$ on
$Y(r')\times \R$, such that the sequence
$\lambda_{r'}(\A_{r'},\Psi_{r'})$ satisfies the following properties. It
converges smoothly on the compact sets $K_V$ to a finite energy solution
$(\A,\Psi)$ of the four-dimensional Seiberg-Witten equations on
$V\times \R$. On the compact sets 
$K_\nu$ it converges to a finite energy solution $(\A,0)$ of the
abelian ASD equation on $\nu(K)\times \R$, and on the compact sets
$K_0$ it converges to a finite energy solution of the
abelian ASD equation on $T^2\times \R^2$.
\end{Pro}

\proof
We can assume that $(\A_r,\Psi_r)$ is in $L^2_{k,\delta_r}(Y(r')\times
\R)$, with $\delta_r\sim 1/r$, by the exponential decay to the
asymptotic values.   On a fixed compact set $K$ the uniform pointwise
bound on the spinor implies the $L^2$-bound
$$ \| \Psi_r \|_{L^2(K)}\leq  Vol(K) s_0. $$
The uniform bound on the curvature provides a uniform $L^2_1$-bound
on the connection (up to a gauge transformation), as in Lemma 5.3.1 of
\cite{Mor}. The uniform pointwise bound on the spinor, together with
the uniform bound
$$ \| \nabla_{\A_r} \Psi_r \|_{L^2(Y(r)\times [t_0,t_1])} \leq
\frac{1}{4} s_0 \| \Psi_r \|_{L^2(Y(r_0)\times [t_0,t_1])} $$
provides a uniform $L^2$-bound for $dF^+_{\A_r}$, as in Lemma 5.3.3 of
\cite{Mor}. A bootstrapping argument then bounds the higher Sobolev
norms. 

\endproof

\subsection{Asymptotics of finite energy ASD connections}

In this subsection we study explicitly the solutions of the abelian
ASD equation on $T^2\times \R^2$ and on an asymptotic end of the form 
$T^2\times [0,\infty)\times \R$.
The first case will provide the geometric limit on compact sets on the
cylinder $T^2\times [-r,r]\times \R$, as $r\to\infty$, and the second
case will give the geometric limit on compact sets on $\nu(K)\times
\R$, and will serve as a model for the asymptotics of the limits on
$V\times \R$, as we discuss later in this work.

If we represent the connection as
$$\A=a(w,s,t)+f(w,s,t)ds +h(w,s,t)dt,$$ 
the abelian ASD equation can be
written as the pair of equations
$$ \partial_t a -dh + *(\partial_s a -df)=0 $$
$$ \partial_t f -\partial_s h + * F_a =0. $$

With a change of variables $z=s+it$ and $z=e^{\rho +i\theta}$, we can
write
$$ a(w,\rho,\theta)=a(w,e^{\rho +i\theta}) $$
$$ f(w,\rho,\theta)=e^{-\rho} \cos\theta \ h(w,e^{\rho
+i\theta})-e^{-\rho}\sin\theta \ f(w,e^{\rho +i\theta}) $$
$$ h(w,\rho,\theta)=e^{-\rho}\cos\theta \ f(w,e^{\rho+i\theta})
+e^{-\rho}\sin\theta \ h(w,e^{\rho+i\theta}). $$

The equations become
$$ \partial_\rho a -d h -*(\partial_\theta a -df)=0 $$
$$ \partial_\rho f -\partial_\theta h -e^{2\rho} * F_a=0. $$
Up to a gauge transformation we can assume that the equations are in
radial gauge, that is $h\equiv 0$ for $\rho$ large enough. Thus, we get
\ba
\partial_\rho a -*(\partial_\theta a -df)=0 
\label{rad1}
\na
\ba
\partial_\rho f -e^{2\rho} * F_a=0. \label{rad2}
\na

Notice that the equation in this form can be interpreted as the
abelian ASD equation on $T^2\times S^1 \times \R$, where the metric on
$T^2$ has a conformal factor depending on $\rho\in \R$, of the form
$e^{-\rho} g_{T^2}$, where $g_{T^2}$ is the standard metric on the
flat torus. In fact the $*$ operator on $p$-forms rescales like
$*_\rho = (e^{-\rho})^{2-2p} *$.

We are considering finite energy solutions, that is, we impose the
condition 
$$ \int_{-\infty}^\infty \int_0^{2\pi} \left( \| \partial_\rho a 
\|^2_{L^2(T^2), g_{T^2}} + e^{2\rho} \| F_a \|^2_{L^2(T^2), g_{T^2}}
\right) d\theta d\rho < \infty. $$

Under a change of variables $\alpha = e^\rho a$ we rewrite equations
(\ref{rad1}) and (\ref{rad2}) as
\ba
\partial_\rho \alpha=\alpha +*\partial_\theta \alpha -e^\rho *df
\label{radt1}
\na
\ba
\label{radt2}
\partial_\rho f= e^{\rho} *d \alpha.
\na

We can write $\alpha = u(w,\rho,\theta) dx + v(w,\rho,\theta) dy$,
where $w=(x,y)$ are the coordinates on the torus $T^2$. We can expand
$u$, $v$, and $f$ in Fourier series in the variables $\theta$, $x$,
and $y$. 

We get
$$ u(w,\theta,\rho)= \sum u_n(\rho)_{lk} \frac{e^{in\theta}}{(2\pi)^{1/2}}
\frac{e^{ilx}}{(2\pi)^{1/2}} \frac{e^{iky}}{(2\pi)^{1/2}}, $$
$$ v(w,\theta,\rho)= \sum v_n(\rho)_{lk} \frac{e^{in\theta}}{(2\pi)^{1/2}}
\frac{e^{ilx}}{(2\pi)^{1/2}} \frac{e^{iky}}{(2\pi)^{1/2}}, $$
$$ f(w,\theta,\rho)= \sum f_n(\rho)_{lk} \frac{e^{in\theta}}{(2\pi)^{1/2}}
\frac{e^{ilx}}{(2\pi)^{1/2}} \frac{e^{iky}}{(2\pi)^{1/2}}. $$

The system of equations becomes the ODE
\ba \frac{d}{d\rho}\left(\begin{array}{c} u_n(\rho)_{lk} \\v_n(\rho)_{lk}
\\ f_n(\rho)_{lk} 
\end{array}\right) =\left(\begin{array}{ccc} 
1 & -in & ike^\rho \\
in & 1 & -il e^\rho \\
ik e^\rho & -il e^\rho & 0 \end{array}\right) \left(\begin{array}{c}
u_n(\rho)_{lk} \\v_n(\rho)_{lk} 
\\ f_n(\rho)_{lk} 
\end{array}\right). \label{ODE:F} \na

In the case of the asymptotics on the end $T^2\times 
[0,\infty)\times \R$, we consider the same
equations and restrict solutions to the domain $\theta \in
[-\pi/2,\pi/2]$, 
in fact, we are not imposing any boundary conditions at $\theta=-\pi/2$ and
$\theta=\pi/2$, other than the functions being smooth across
$\theta=-\pi/2$ and $\theta=\pi/2$.

The assumption that $\A$ is a $U(1)$ connection imposes the constraint
on the coefficients
$$ u_n(\rho)_{lk}=-\overline {u_{-n}(\rho)_{-l,-k}}, $$
$$ v_n(\rho)_{lk}=-\overline { v_{-n}(\rho)_{-l,-k}}, $$
$$ f_n(\rho)_{lk}=-\overline {f_{-n}(\rho)_{-l,-k}}. $$

A direct analysis of this system of ODE's (cf. \cite{Ha}, \S X) proves the
following Proposition.

\begin{Pro}
The only finite energy solutions of the abelian ASD equation on
$T^2\times \R^2$ are flat connections on $T^2$, constant in the $\R^2$ 
directions. On the asymptotic end $T^2\times
[r_0,\infty)\times \R$  all the finite energy
solutions are of the form 
\ba
\label{allfinenergy}
\begin{array}{lr}
a(w,\rho,\theta)=& ( u_0 + \sum_{n\neq 0} u_n e^{-|n|\rho+in\theta}  \\[2mm]
& +\sum_{n,l,k} lc_{nlk} e^{i(lx+ky+n\theta)}) dx +\\[2mm]
& ( v_0 + \sum_{n\neq 0} v_n e^{-|n|\rho+in\theta}  \\[2mm]
& +\sum_{n,l,k} kc_{nlk} e^{i(lx+ky+n\theta)} )dy \\[4mm]
f(w,\rho,\theta)=&  f_0 + \sum_{n,l,k} n c_{nlk} e^{i(lx+ky+n\theta)} 
\end{array}
\na
with the coefficients that satisfy
$$ \sum_n |n u_n|^2 + \sum_{n,l,k} |(n^2+l^2+k^2) c_{nlk} |^2 <\infty, $$
and
$$ c_{nlk}=-\overline{c_{-n,-l,-k}} \ \ \hbox{and} \ \
u_n=-\overline{u_{-n}}. $$

In radial gauge, the limit in the $\rho\to\infty$ direction is given
by 
$$ \begin{array}{lr}
a_\infty(w,\theta)= & (u_0 + \sum_{n,l,k} lc_{nlk}
e^{i(lx+ky+n\theta)})dx + \\[2mm]
& (v_0 + \sum_{n,l,k} kc_{nlk} e^{i(lx+ky+n\theta)})dy  \\[4mm]
f_\infty(w,\theta) = & f_0 + \sum_{n,l,k} n c_{nlk} e^{i(lx+ky+n\theta)}.
\end{array} $$
\end{Pro}

\proof

It is easy to find a special family of solutions of the form
$$ c_{nlk} \left(\begin{array}{c} l\ e^\rho \\ k e^\rho \\ n
\end{array}\right), $$ 
with $-\overline{c_{-n,-l,-k}}=c_{nlk}$.
These correspond to the condition $d\alpha =0$ in the equation
(\ref{radt2}). In the original variables, that is multiplying the
connection terms by $e^{-\rho}$, these are solutions of (\ref{rad1})
and (\ref{rad2}) constant in the 
$\rho$-direction. 

We assume for the moment that $(k,l)\neq (0,0)$. Following the
standard methods of ODE theory, we use these solutions in order to
reduce the system to a system of two equations. If we assume $l\neq
0$, we can consider the matrix
$$ Z=\left( \begin{array}{ccc} le^\rho & 0&0\\
ke^\rho & 1&0\\
n &0&1 \end{array}\right). $$

The change of variables
$$ \left( \begin{array}{c} u\\v\\f\end{array}\right)=Z\left(
\begin{array}{c} \omega\\ \nu\\ \phi \end{array}\right) $$
gives the new system of equations
$$ \frac{d}{d\rho}\left( \begin{array}{c}\omega\\ \nu\\ \phi
\end{array}\right) = Z^{-1} A \left( \begin{array}{ccc} 0&0&0\\0&1&0\\
0&0&1 \end{array}\right)\left(\begin{array}{c} \omega\\ \nu\\ \phi
\end{array}\right), $$
where $A$ is the matrix
$$ A=\left(\begin{array}{ccc} 
1 & -in & 0 \\
in & 1 & 0 \\
0 & 0 & 0 \end{array}\right). $$ 

This system consists of the equation
$$ \omega' = \frac{-in}{l}e^{-\rho} \nu + \frac{ik}{l} \phi $$
and the two by two system
$$ \frac{d}{d\rho} \left( \begin{array}{c} l\nu\\ l\phi
\end{array}\right)=\left( \begin{array}{cc} ink+l & -i(k^2+l^2)e^\rho 
\\ i(n^2 e^{-\rho} -l^2 e^\rho) & -ink \end{array}\right) \left(
\begin{array}{c} \nu\\ \phi \end{array}\right). $$

Since we are interested in the large $\rho >> 0$ behavior of the
system we can isolate a leading term and treat the rest of the system
as a perturbation.

Let us define 
$$ L_\rho= \left( \begin{array}{cc} 0 & -i(k^2+l^2)e^\rho \\ -i 
l^2 e^\rho & 0 \end{array}\right) $$
to be the leading term and
$$ P_\rho= \left( \begin{array}{cc} ink+l & 0 \\
in^2 e^{-\rho} & -ink \end{array}\right) $$
to be the perturbation.

The unperturbed system has eigenvalues 
$$ \lambda_\pm (\rho) = \pm i e^\rho (k^2+l^2)^{1/2}, $$
with eigenvectors
$$ U_\pm = \left( \begin{array}{c} 1 \\ \frac{\pm l}{(k^2+l^2)^{1/2}}
\end{array}\right). $$
So, upon diagonalizing the matrix we obtain solutions
$\exp(\pm ie^\rho (k^2+l^2)^{1/2})$, hence in the original system of
coordinates we have solutions
\ba \begin{array}{rl}
\nu(\rho)= & \frac{(k^2+l^2)^{1/2}}{2l} c_1 \exp(ie^\rho
(k^2+l^2)^{1/2}) \\[2mm] & + \frac{(k^2+l^2)^{1/2}}{2l} c_2
\exp(-ie^\rho (k^2+l^2)^{1/2}),  \end{array} 
\label{2x2,1}
\na
and
\ba \begin{array}{rl}
\phi(\rho)= & -\frac{1}{2} c_1 \exp(ie^\rho
(k^2+l^2)^{1/2}) \\[2mm]
& + \frac{1}{2} c_2 \exp(-ie^\rho (k^2+l^2)^{1/2}). \end{array} 
\label{2x2,2}
\na
with the remaining equation that gives
\ba
\begin{array}{l}
\omega(\rho) =\int_0^\rho \left( \frac{-in(k^2+l^2)^{1/2}}{l^2} (
c_1 \exp(ie^\tau) + c_2\exp(-ie^\tau) )\right. - \\[2mm]
\left. \frac{ik}{2l}e^\tau (c_1
\exp(ie^\tau) - c_2\exp(-ie^\tau) ) \right) d\tau. \end{array}
\label{2x2,3}
\na

The perturbed system has eigenvalues 
$$ \tilde\lambda_\pm (\rho) = \frac{1\pm \sqrt{1-4(e^{2\rho}(k^2+l^2) -n^2
-\frac{ink}{l})} }{2}. $$
The matrix can be diagonalized so that the solutions are of the form
$$ \exp\left(\int_0^\rho \tilde\lambda_\pm(\tau) d\tau\right).
$$
The long distance $\rho>>0$ behavior of these solutions is given by
the asymptotics $e^{\rho \pm ie^\rho}$.
On the other hand the eigenvectors become asymptotically
$\rho$-independent and approach the eigenvectors $U_\pm$ of the unperturbed
system. Thus, the asymptotic behavior of the solutions will be of the 
form
\ba
\label{2x2,sol1}
\nu(\rho)\sim \frac{(k^2+l^2)^{1/2}}{2l} \left( c_1 e^{\rho + ie^{\rho}}
+ c_2 e^{\rho - ie^{\rho}} \right)
\na
and
\ba
\label{2x2,sol2}
\phi(\rho)\sim \frac{-1}{2} \left( c_1 e^{\rho + ie^{\rho}}
- c_2 e^{\rho - ie^{\rho}} \right)
\na
The third variable is then obtained as
\ba
\label{2x2,sol3}
\begin{array}{l}
\omega(\rho) \sim\int_0^\rho \left( \frac{-in(k^2+l^2)^{1/2}}{l^2} (
c_1 \exp(ie^\tau) + c_2\exp(-ie^\tau) )\right. - \\[2mm]
\left. \frac{ik}{2l}e^\tau (c_1
\exp(ie^\tau) - c_2\exp(-ie^\tau) ) \right) d\tau. \end{array}
\na
Then the original solutions will be of the form
$$ \left(\begin{array}{c} u \\ v \\ f
\end{array}\right)=\left(\begin{array}{c} l e^\rho \omega(\rho) \\
ke^\rho \omega(\rho)+\nu(\rho) \\ n \omega(\rho) +\phi(\rho)
\end{array}\right). $$

The case with $l=0$ and $k\neq 0$ is analogous. In fact, in that case
we can use a similar reduction by considering the matrix
$$ Z=\left( \begin{array}{ccc} 0 & 1&0\\
ke^\rho & 0&0\\
n &0&1 \end{array}\right). $$
The change of variables
$$ \left( \begin{array}{c} u\\v\\f\end{array}\right)=Z\left(
\begin{array}{c} \omega\\ \nu\\ \phi \end{array}\right) $$
gives the new system of equations
$$ \frac{d}{d\rho}\left( \begin{array}{c}\omega\\ \nu\\ \phi
\end{array}\right) = Z^{-1} A \left( \begin{array}{ccc} 0&1&0\\0&0&0\\
0&0&1 \end{array}\right)\left(\begin{array}{c} \omega\\ \nu\\ \phi
\end{array}\right), $$ 
which is of the form
$$ \omega' = \frac{in}{k}e^{-\rho} \nu  $$
and the remaining two by two system
$$ \frac{d}{d\rho} \left( \begin{array}{c} \nu\\ \phi
\end{array}\right)=\left( \begin{array}{cc} 1 & ike^\rho 
\\ \frac{-in^2}{k} e^{-\rho} +ik e^\rho & 0 \end{array}\right) \left(
\begin{array}{c} \nu\\ \phi \end{array}\right). $$
Again we can isolate a leading term
$$ L_\rho=\left( \begin{array}{cc} 0 & ike^\rho \\ ik e^\rho & 0
\end{array}\right) $$
and a perturbation
$$ P_\rho=\left( \begin{array}{cc}  1 & 0 \\ \frac{-in^2}{k}
e^{-\rho} & 0 \end{array}\right). $$
The unperturbed system has eigenvalues $\lambda_\pm= \pm ike^\rho$ and 
$\rho$-independent eigenvectors 
$$ U_\pm = \left( \begin{array}{c} 1 \\ \pm 1 \end{array}\right). $$
The perturbed system has eigenvalues 
$$\tilde \lambda_\pm =\frac{1\pm\sqrt{1-4(k^2e^{2\rho}-n^2)}}{2}. $$
The asymptotics of the solutions is $\exp(\rho \pm ike^{\rho})$, with
eigenvectors that asymptotically approach the $\rho$-independent
eigenvectors $U_\pm$.
Thus, we obtain asymptotics 
$$ \omega(\rho)\sim \frac{in}{2k} \int_0^\rho (c_1 e^{+ ike^{\tau}} +
c_2 e^{- ike^{\tau}}) d\tau, $$
$$ \nu(\rho)\sim \frac{1}{2} (c_1 e^{\rho + ike^{\rho}} + c_2 e^{\rho
- ike^{\rho}} )$$
and
$$ \phi(\rho)\sim \frac{1}{2} (c_1 e^{\rho + ike^{\rho}} - c_2 e^{\rho
- ike^{\rho}} ).$$

Among these families of solutions, the only elements that satisfy the
finite energy condition are written in the original variables as
solutions of (\ref{rad1}) and
(\ref{rad2}) of the form
\ba \label{onefinenergy}
\begin{array}{lr}
a(w,\rho,\theta)=& ( u_0 
 +\sum_{n,l,k} lc_{nlk} e^{i(lx+ky+n\theta)}) dx +\\
& ( v_0  +\sum_{n,l,k} kc_{nlk} e^{i(lx+ky+n\theta)} )dy \\[3mm]
f(w,\rho,\theta)=&  f_0 + \sum_{n,l,k} n c_{nlk} e^{i(lx+ky+n\theta)}. 
\end{array}
\na

The remaining solutions satisfy $(k,l)=(0,0)$. These satisfy
the conditions $df=0$ and $d\alpha=0$ in the equations (\ref{radt1})
and (\ref{radt2}).
These are solutions of the form
$$ \left(\begin{array}{c} c_1(n,0,0) e^{(1+n)t} - ic_2(n,0,0) e^{(1-n)t} \\
c_1(n,0,0) e^{(1+n)t} + ic_2(n,0,0) e^{(1-n)t} \\
c_3(n,0,0) \end{array}\right) $$
to the system with constant coefficients with matrix
$$ A=\left(\begin{array}{ccc} 
1 & -in & 0 \\
in & 1 & 0 \\
0 & 0 & 0 \end{array}\right). $$

This gives solutions of the original equations (\ref{rad1}) and
(\ref{rad2}) of the form 
\ba
\label{finen}
\begin{array}{ll}
a(w,\rho,\theta)=& \left( u_0 + \sum_{n\neq 0} u_n e^{-|n|
\rho+in\theta} \right) dx  \\ & + \left( v_0+ \sum_{n\neq 0} v_n e^{-|n|
\rho+in\theta} \right) dy \\
f(w,\rho,\theta)= & f_0 \end{array}
\na

This proves the Proposition. In fact, we can write the function 
$$f_\infty -f_0 = -i\frac{d}{d\theta} \sum_{n,l,k} c_{nlk}
e^{i(lx+ky+n\theta)}, $$
so that, if we define 
$$ \gamma(\theta)=\sum_{l,k} -ic_{nlk} e^{i(lx+ky+n\theta)}, $$
the flat connection $a_\infty$ can be written as
$$ a_\infty - (u_0 dx + v_0 dy) = d \gamma(\theta). $$
This means that the path of asymptotic flat connections in the covering
$\chi_0(T^2,\nu(K))$ of the character variety $\chi(T^2)$ can be written as
$$ a_\infty (\theta) =\lambda(\theta) \cdot a_0, $$
where $\lambda: T^2\times [0,\pi]\to U(1)$ satisfies
$$\lambda^{-1}(\theta)d_{T^2}\lambda(\theta)=d_{T^2}\gamma(\theta). $$ 

\endproof

Thus we have obtained the following.

\begin{Cor}
A finite energy solution of the abelian ASD equation on an asymptotic
end of the form $T^2\times [0,\infty)\times \R$ decays
exponentially along the radial direction. The asymptotic values for
different angles $\theta$ are all
within the same gauge class of flat connections on $T^2$. Thus, we
have decay to a point $a_\infty = [a_\infty (\theta)]$.
\end{Cor}

\begin{Rem}
Up to a global change of gauge, the only solutions of the form
(\ref{allfinenergy}) that extend to all of $\nu(K)\times \R$ are flat
connections on $T^2$ constant in the $s$ and $t$ directions.   
\label{constant}
\end{Rem}

In the light of Remark \ref{constant}, the explicit analysis of the
asymptotics (\ref{allfinenergy}) seems a somewhat useless
complication. However, the analysis above 
will be crucial in the following section, in order to describe, with a
perturbative analysis, the asymptotic
behavior of finite energy solutions of the Seiberg-Witten equations
on $V\times \R$ with an end $T^2\times
[0,\infty)\times \R$.

\subsection{Asymptotics of monopoles on $V\times \R$}

We now analyze the asymptotics of solutions of the four-dimensional
Seiberg-Witten equations on the end $T^2\times [r_0,\infty)\times \R$
of the manifold $V\times \R$, where $V$ is the knot complement in $V$,
endowed with an infinite cylindrical end.

We write the connection $\A=a(w,s,t)+f(w,s,t)ds +h(w,s,t)dt$, and the
spinor section $\Psi=(\alpha,\beta)$ with $\alpha(s,t)\in
\Lambda^{0,0}(T^2)$ and $\beta(s,t)\in \Lambda^{0,1}(T^2)$.

The Seiberg-Witten equations can be written in the form
$$ \partial_t a -dh + *(\partial_s a -df)= * i(\bar\alpha \beta
+\alpha\bar\beta) $$ 
$$ \partial_t f -\partial_s h + * F_a
=\frac{i}{2}(|\alpha|^2-|\beta|^2), $$
for the curvature equation, and 
$$ \partial_t \alpha + h\alpha + i \partial_s \alpha + if \alpha +
\bar\partial_a^*\beta =0 $$
$$ \partial_t \beta + h\beta -i \partial_s \beta -if\beta +
\bar\partial_a \alpha =0, $$
for the Dirac equation.

As before, we introduce the variables $w=(x,y)\in T^2$, and $z=s+it$,
$z=e^{\rho +i\theta}$ on $[r_0,\infty)\times \R$, and the change of
coordinates 
\ba \begin{array}{c}
 a(w,\rho,\theta)=a(w,e^{\rho +i\theta}) \\[2mm]
f(w,\rho,\theta)=e^{-\rho} \cos\theta \ h(w,e^{\rho
+i\theta})-e^{-\rho}\sin\theta \ f(w,e^{\rho +i\theta}) \\[2mm]
h(w,\rho,\theta)=e^{-\rho}\cos\theta \ f(w,e^{\rho+i\theta})
+e^{-\rho}\sin\theta \ h(w,e^{\rho+i\theta}) \\[2mm]
 \alpha(w,\rho,\theta)= \alpha(w,e^{\rho +i\theta}) \\[2mm]
 \beta(w,\rho,\theta)=\beta(w,e^{\rho +i\theta}). \end{array}
\label{change:polar} \na

The Seiberg-Witten equations in radial gauge (i.e. with $h\equiv 0$
for large $\rho$) are then written in the form
\ba
\label{SWradial}
\begin{array}{l}
\partial_\rho a =*(\partial_\theta a -df + i(\bar\alpha \beta
+\alpha\bar\beta)) \\[2mm]
\partial_\rho f =e^{2\rho} *( F_a +
\frac{i}{2}(|\alpha|^2-|\beta|^2)\omega) \\[2mm]
\partial_\rho \alpha = i( \partial_\theta \alpha + f \alpha + e^{\rho
+i\theta} \bar\partial_a^*\beta) \\[2mm]
\partial_\rho \beta = -i( \partial_\theta \beta + f\beta +
e^{\rho-i\theta} \bar\partial_a \alpha). \end{array}
\na

In order to study the asymptotic behavior of these solutions, we can
again use Fourier transform. We first analyze the asymptotic behavior 
of a linear system and then introduce the non-linear terms in 
a sequence of successive approximations \cite{Ha}.
In order to simplify the expression of the quadratic terms in
(\ref{SWradial}), we use the notation
$$ A\cdot B =\sum (a b)_{nlk} \frac{e^{i(n\theta + lx
+ky)}}{(2\pi)^{3/2}}, $$
with the coefficients given by
$$ (a b)_{nlk}=\sum \frac{1}{(2\pi)^{3/2}} a_{ijh} b_{n-i \ l-j \ k-h}, $$  
for
$$ A=\sum a_{nlk} \frac{e^{i(n\theta + lx
+ky)}}{(2\pi)^{3/2}} $$
and
$$ B=\sum b_{nlk} \frac{e^{i(n\theta + lx
+ky)}}{(2\pi)^{3/2}}. $$

With this notation, the equation (\ref{SWradial}) becomes
\ba
\label{SWradfourier}\begin{array}{ll}
u_{nlk}'=& -in v_{nlk} + ik  f_{nlk} + 2i Re
(\bar\alpha \beta)_{nlk} \\[2mm]
v_{nlk}'=& in u_{nlk}  -il f_{nlk} +2i
Im(\bar\alpha\beta)_{nlk} \\[2mm] 
f_{nlk}'=& ik e^{2\rho} u_{nlk} -il e^{2\rho} v_{nlk} + e^{2\rho} (\bar
\alpha \alpha)_{nlk} - (\bar\beta\beta)_{nlk} \\[2mm]
\alpha_{nlk}'=&  -n\alpha_{nlk} + i e^\rho \frac{1}{2}(il-k) \bar\beta_{-n-1 
\ -l \ -k}+ \\[2mm]
& i e^\rho \frac{1}{2}(a\bar\beta)_{nlk} + i
(f\alpha)_{nlk}\\[2mm] 
\beta_{nlk}'=& n \beta_{nlk} -i e^\rho (il-k) \alpha_{n+1 \ lk} \\[2mm]
& -ie^\rho (a\alpha)_{nlk} -i (f\beta)_{nlk}. \end{array}
\na

Consider the linear system of equations
\ba
\label{SWradlin}\begin{array}{l}
\partial_\rho a =*(\partial_\theta a -df)\\[2mm]
\partial_\rho f =e^{2\rho} * F_a\\[2mm]
\partial_\rho \alpha = i( \partial_\theta \alpha + e^{\rho
+i\theta} \bar\partial^*\beta) \\[2mm]
\partial_\rho \beta = -i( \partial_\theta \beta + 
e^{\rho-i\theta} \bar\partial \alpha). \end{array}
\na

This is the uncoupled system of the abelian ASD equation and the
linear Dirac equation. We analyzed the first two equations 
in the previous section. 

In the Dirac equations, since again we are interested in the large
$\rho >> 0$ behavior of the finite energy solutions, we can isolate a
leading term and a perturbation. 

The leading term gives a system of the form
\ba \begin{array}{l} 
\alpha_{nlk}'= i e^\rho \frac{1}{2}(il-k) \bar\beta_{-n-1 
\ -l \ -k}\\[2mm]
\beta_{nlk}'=-i e^\rho (il-k) \alpha_{n+1 \ lk} , \end{array}
\label{leading:spinor} \na
Whereas we identify perturbation all the terms of the form
\ba \left( \begin{array}{cc} -n & 0 \\[2mm] 0 & n \end{array} \right)
\left(\begin{array}{c} \alpha_{nlk} \\[2mm] \beta_{nlk}\end{array}
\right). \label{pert:spinor} \na

We analyze solutions to the unperturbed system (\ref{leading:spinor}).
If we introduce the notation $\alpha_{nlk}=\eta_{nlk}+ i \xi_{nlk}$ and
$\beta_{nlk}= p_{nlk} + i q_{nlk}$, the system (\ref{leading:spinor})
uncouples in the independent systems of eight equations: 
\ba \label{uncoupled}\begin{array}{l}
p_{nlk}'= e^\rho( l \eta_{n+1 \ lk} - k \xi_{n+1 \ lk}) \\[2mm]
q_{nlk}'= e^\rho( l \xi_{n+1 \ lk} + k \eta_{n+1 \ lk}) \\[2mm]
\eta_{n+1 \ lk}'= \frac{e^\rho}{2}( -l p_{-n-2 \ -l \ -k} -k q_ {-n-2
\ -l \ -k}) \\[2mm]
\xi_{n+1 \ lk}'= \frac{e^\rho}{2}( l p_{-n-2 \ -l \ -k} -k q_ {-n-2
\ -l \ -k}) \\[2mm]
p_{-n-2 \ -l \ -k}'= e^\rho( -l \eta_{-n-1 \ -l \ -k} +k \xi_{-n-1 \
-l \ -k}) \\[2mm]
q_{-n-2 \ -l \ -k}'= e^\rho( -l \eta_{-n-1 \ -l \ -k} -k \xi_{-n-1 \
-l \ -k}) \\[2mm]
\eta_{-n-1 \ -l \ -k}'=\frac{e^\rho}{2}( l p_{nlk} +k q_{nlk}) \\[2mm]
\xi_{-n-1 \ -l \ -k}'=\frac{e^\rho}{2}( -l p_{nlk} +k q_{nlk}) 
\end{array} \na

Consider the matrices
$$ A=\left( \begin{array}{cc} l& -k\\ l&k \end{array}\right) $$
and
$$ B =\frac{1}{2} \left( \begin{array}{cc} -l & -k \\ l & -k
\end{array}\right).  $$

Up to a reparametrization in the variable $\tau= e^\rho$, we obtain
the autonomous linear system with the $8\times 8$ matrix

$$ M=\left( \begin{array}{cccc} 0&A&0&0 \\ 0&0&B&0 \\ 0&0&0&-A \\
-B&0&0&0 \end{array} \right). $$

If we write $\eta_0=(k^2 +6kl + l^2)^{1/2}$, $\eta_1=k-l$, and
$\eta_2=k+l$, then the eigenvalues of $M$ are of the form $\pm
\lambda^{nlk}_i$, for $i=1,\ldots 4$, with
\ba \begin{array}{lr}
\lambda^{nlk}_1=\frac{1}{2} (\eta_2^2 +\eta_1\eta_0)^{1/2} &
\lambda^{nlk}_2=\frac{1}{2} (\eta_2^2 -\eta_1\eta_0)^{1/2}  \\[2mm]
\lambda^{nlk}_3=\frac{1}{2} (-\eta_2^2 -\eta_1\eta_0)^{1/2} &
\lambda^{nlk}_4=\frac{1}{2} (-\eta_2^2 +\eta_1\eta_0)^{1/2}.
\end{array} \label{eigenv} \na

It is easy to see that, in the basis $U_i^\pm(nlk)$ of eigenvectors,
the solutions corresponding to $\lambda^{nlk}_1$ and $\lambda^{nlk}_2$ 
are saddles for all non-trivial $(n,l,k)$ and solutions corresponding
to $\lambda^{nlk}_3$ and $\lambda^{nlk}_4$ are centers for all
non-trivial $(n,l,k)$. Thus, the finite energy solutions will be of
the form
\ba \left(\begin{array}{c}
p_{nlk}\\ q_{nlk}\\ \eta_{n+1 \ lk}\\ \xi_{n+1 \ lk}\\ p_{-n-2 \ -l \
-k} \\ q_{-n-2 \ -l \ -k} \\ \eta_{-n-1 \ -l \ -k}\\ \xi_{-n-1 \ -l \ -k}
\end{array}\right) = \sum_{i=1}^2 c_i^{nlk} \exp(-|\lambda_i^{nlk}|\tau)
U_i^-(nlk) , \label{U} \na
for $\tau=e^{\rho}$.

If we add the perturbation terms (\ref{pert:spinor}) we can proceed as
in the case of the ASD equation and observe that, for large
$\rho\to\infty$ the eigenvalues and eigenvectors of the perturbed
system converge to the eigenvalues and eigenvectors
$\pm\lambda^{nlk}_i$ and $U_i^\pm(nlk)$ described above. 

In the variables
$\alpha_{nlk}=\eta_{nlk}+ i \xi_{nlk}$ and 
$\beta_{nlk}= p_{nlk} + i q_{nlk}$, the perturbed system can be
written as the non-autonomous linear system with the $8\times 8$ matrix
\ba \tilde M=\left( \begin{array}{cccc} n I &e^\rho A&0&0 \\
0&-(n+1)I&e^\rho B&0 \\ 0&0&-(n+2)I&-e^\rho A \\ 
-e^\rho B&0&0&(n+1) I \end{array} \right), \label{tildeM} \na
where $I$ is the $2\times 2$ identity matrix.

The eigenvalues have a term independent of $\rho$ which reduces to
the diagonal entries of the matrix $\tilde M$ as $\rho\to -\infty$, and
terms containing  $e^\rho$, of which the highest order term reduces to
the eigenvalues of the matrix $M$. Thus, the rate of decay as
$\rho\to\infty$ in the directions corresponding to eigenvalues of
negative real part is at least given by $\exp(-n \rho)$ and at
most by $\exp(-C_{nlk} e^\rho)$, for some positive constants
$C_{nlk}>0$, or by mixed terms  $\exp(n\rho -C_{nlk} e^\rho)$.

With this analysis in place, we know the  asymptotic behavior as
$\rho\to\infty$ of the finite energy solutions of the linear system. 
Now we want to analyze the asymptotic behavior of solutions of the
original non-linear equations.

Suppose given a finite energy solution
$\Xi_0=(u_0,v_0,f_0,\alpha_0,\beta_0)$ of the linear system. Consider
a perturbed linear system with matrix $L_{\Xi_0}$ as follows:
\ba
\label{pert0}\begin{array}{ll}
u_{nlk}'=& -in v_{nlk} + ik  f_{nlk} + 2i Re
(\bar\alpha_0\beta)_{nlk} +2iRe(\bar\alpha\beta_0)_{nlk}\\[2mm]
v_{nlk}'=& in u_{nlk}  -il  f_{nlk} +2i
Im(\bar\alpha_0\beta)_{nlk}+2iIm(\bar\alpha\beta_0)_{nlk} \\[2mm] 
f_{nlk}'=& ik e^{2\rho} u_{nlk} -il e^{2\rho} v_{nlk} \\[2mm] & + 2e^{2\rho}(
Re(\bar 
\alpha_0\alpha)_{nlk} - (\bar\beta_0\beta)_{nlk}) \\[2mm]
\alpha_{nlk}'= & -n\alpha_{nlk} + i e^\rho \frac{1}{2}(il-k) \bar\beta_{-n-1 
\ -l \ -k}\\[2mm] & + i e^{\rho} \frac{1}{2}(a_0\bar\beta)_{nlk}+ i
e^{\rho} 
\frac{1}{2}(a\bar\beta_0)_{nlk} + i (f_0\alpha)_{nlk}+
i(f\alpha_0)\\[2mm] 
\beta_{nlk}'=& n \beta_{nlk} -i e^\rho (il-k) \alpha_{n+1 \ lk}
-ie^{\rho} (a_0\alpha)_{nlk} -ie^{\rho} (a\alpha_0)_{nlk} \\[2mm] & -i
(f_0\beta)_{nlk}-i(f\beta_0) . \end{array}
\na

Similarly, we can define inductively solutions $\Xi_{\nu+1}$
satisfying the linear system of equations
$$ \Xi_{\nu+1} ' =L_{\Xi_\nu} \Xi_{\nu+1}, $$
with initial condition $\Xi_{\nu+1}(0)=\Xi_\nu(0)$. 

\begin{Pro}
All the solutions $\Xi_{\nu+1}$ are finite energy. Moreover, we have
$$ \Xi_\nu \to \Xi $$
uniformly on compact sets along with all derivative, where $\Xi$ is a
solution of the original non-linear system (\ref{SWradial})
$$ \Xi ' =L_{\Xi} \Xi. $$ 
The solution $\Xi$ has the same asymptotic behavior as
$\rho\to\infty$ as $\Xi_0$, up to terms that decay faster as
$\rho\to\infty$. 
\end{Pro}

\proof
Inductively, the system of equations satisfied by $\Xi_{\nu+1}$ can be
separated into a leading term and a perturbation. The rate of decay
of the solutions to the system with the leading term only will
determine the asymptotics of $\Xi_{\nu+1}$.

In the case of the system (\ref{pert0}), we have the coefficients of
$\Xi_0$ decaying at least like $e^{-n\rho}$. Thus, for large $\rho$,
the leading terms in the system (\ref{pert0}) will be the same as the
leading terms in the original linear system.

Inductively, at the $\nu$-th stage, we have the original linear system 

\ba
\begin{array}{ll}
u_{nlk}'=& -in v_{nlk} + ik  f_{nlk} \\ & 2i Re
(\bar\alpha_\nu\beta)_{nlk} +2iRe(\bar\alpha\beta_\nu)_{nlk}\\[2mm]
v_{nlk}'=& in u_{nlk}  -il  f_{nlk} \\ & 2i
Im(\bar\alpha_\nu\beta)_{nlk}+2iIm(\bar\alpha\beta_\nu)_{nlk} \\[2mm] 
f_{nlk}'=& ik e^{2\rho} u_{nlk} -il e^{2\rho} v_{nlk} \\ & 2e^{2\rho}(
Re(\bar \alpha_\nu\alpha)_{nlk} - (\bar\beta_\nu\beta)_{nlk}) \\[2mm]
\alpha_{nlk}'= & -n\alpha_{nlk} + i e^\rho \frac{1}{2}(il-k) \bar\beta_{-n-1 
\ -l \ -k}\\ & i e^{\rho} \frac{1}{2}(a_\nu\bar\beta)_{nlk}+ i
e^{\rho} 
\frac{1}{2}(a\bar\beta_\nu)_{nlk} + i (f_\nu\alpha)_{nlk}+
i(f\alpha_\nu)\\[2mm] 
\beta_{nlk}'=& n \beta_{nlk} -i e^\rho (il-k) \alpha_{n+1 \ lk} \\ &
-ie^{\rho} (a_\nu\alpha)_{nlk} -ie^{\rho} 
(a\alpha_\nu)_{nlk} -i 
(f_\nu\beta)_{nlk}-i(f\beta_\nu) . \end{array}
\label{pert0I}
\na

Under the inductive hypothesis that the coefficients of the solutions
$\Xi_\nu$ decay at a rate at least $e^{-n\rho}$ as $\rho\to\infty$, we
see that the leading terms in the system (\ref{pert0I}) are the same
as in the original unperturbed system. The eigenvalues
$\pm\tilde\lambda_i^{nkl}(\rho)$ of the perturbed 
system (\ref{pert0I}) and the corresponding eigenvectors $\tilde
U_i^\pm(\rho)$ converge asymptotically to the $\rho$-independent
eigenvalues and eigenvectors of the matrix $M$. Thus, the solutions
$\Xi_\nu$ are of finite energy and 
maintain the same rate of decay as the original solution $\Xi_0$.

\endproof

\subsection{Uniform convergence}

Notice that, in the particular case where the pointwise convergence
proved in Lemma
\ref{pointwise:conv} extends uniformly to an infinite strip $T^2\times 
[-r_0,r_0]\times \R$, in the sense specified below, then we have
uniform convergence of the $(\A_r,\Psi_r)$ in the smooth topology on
the infinite strip to a flat connection. In particular, if we consider 
the geometric limits of solutions on (\ref{3SW}) on $Y(r)$ as
$r\to\infty$, as analysed in \cite{CMW}, we see that, in this special
case, the asymptotic  
values $(A_r(\pm\infty),\psi_r(\pm\infty))$ break through the same
flat connection on $T^2$, when $r\to\infty$.

By the results of the previous section, on any fixed
compact set $K$ in $T^2\times [-r,r]\times \R$, or in $\nu(K)\times \R$, 
we have 
\[ \| F_{\A_r} \|_{L^2(K)} \to 0. \]

\begin{Pro}
Suppose given a family of finite energy solutions
$\Xi_r=(\A_r,\Psi_r)$, for all $r\geq r_0$. 
Suppose there exists a
strip $T^2\times [r_0, r_1 ]\times \R$ in $T^2\times [-r,r]\times \R$
such that we have 
\ba
\| F_{\A_r} \|_{L^2(T^2\times [r_0, r_1 ]\times \R)} \to 0 
\label{uniform:curv} \na
as $r\to\infty$. Then there exists a flat connection $a_{T^2}$ on
$T^2$ such that (up to gauge transformations) the sequence $\Xi_r$
converges in the $L^2_2$-topology on $T^2\times [r_0, r_1 ]\times \R$
to $a_{T^2}$. In particular this implies that,
under the splitting of the 
critical points ${\cal M}_{Y(r)}={\cal M}_{V}\#_{\chi_0(T^2)}
\chi(\nu(K))$, for large enough $r$, the asymptotic values
$(A_r(\pm\infty),\psi_r(\pm\infty))$ can be written as the gluing
$$(A_r(\pm\infty),\psi_r(\pm\infty))=
(A^\prime(\pm\infty),\psi^\prime(\pm\infty)) \#_{a_{T^2}} (\lambda^\pm
a_{T^2}, 0), $$
that is, they split through the same asymptotic flat connection on
$T^2$. 
Moreover, suppose that $(\A',\Psi')$ and $(\A,0)$ are the limits on
compact sets on $V\times \R$ and $\nu(K)\times \R$ respectively, with
asymptotic values $[a_V] \in \chi_0(T^2,V)$ and $[a_\nu] \in
\chi_0(T^2,\nu(K))$, as described in the previous sections. Then, in
this particular case, there 
are gauge transformations $\lambda'$ and $\lambda$ on $V$ and $\nu(K)$ 
respectively such that we have
$$ \lambda' a_V= a_{T^2} = \lambda a_\nu. $$
\label{unif:conv}
\end{Pro}

\proof
By Uhlenbeck weak convergence \cite{FU} \cite{U}, and the result of
Proposition 8.3 of \cite{FU}, the condition
$$ \| F_{\A_r} \|_{L^2(T^2\times [r_0, r_1 ]\times \R)} \to 0  $$
implies that there exists a flat connection
$\A$ on $T^2\times [r_0, r_1 ]\times \R$ such that, up to gauge
transformations $\lambda_r$ in 
$L^2_3$ and up to passing to a subsequence, we have
\[ \| \lambda_r\A_r -\A \|_{L^2_2 (T^2\times [r_0, r_1 ]\times
\R)}\leq C \| F_{\A_r} \|_{L^2(T^2\times [r_0, r_1 ]\times
\R)} \to 0. \]  
We can now show that we can further gauge transform the connections
$\A_k$ so that this Uhlenbeck limit on $ T^2\times [r_0, r_1 ]\times
\R$ is a flat connection on $T^2$, constant in the $t\in \R$ and $s\in
[r_0,r_1]$ directions. 

{\bf Claim:} Let $\A$ be an $L^2_{2,loc}$ flat connection on
$T^2\times [r_0, r_1 ]\times\R$. Then, 
there exist and $L^2_{3,loc}$ gauge transformation $\lambda$ on
$T^2\times [r_0, r_1 ]\times \R$ such that $\lambda\A$ is a flat
connection $a_{T^2}$ on 
$T^2$, constant in the $t$ and $s$ directions. 

{\bf proof of Claim:} We can write $\A$ in the
form $\A=a_{T^2} + if\ dt +ih\ ds$, where $a_{T^2}$ is a connection on
$T^2$, and $f(x,y,t,s)$ and $h(x,y,t,s)$ are real valued functions.
Then the equation $F_{\A}=0$ becomes
\ba
\begin{array}{c} F_{a_T^2}=0 \\ \partial_t a_{T^2}=d_{T^2} if \\
\partial_s 
a_{T^2}= d_{T^2} ih \\ \partial_t h= \partial_s f. \end{array}
\label{gaugeF}  \na
 
Define the gauge transformation
$$ \lambda = \exp \left( -i \int_{r_0}^s h(x,y,t,\sigma)d\sigma -i
\int_0^t f(x,y,\tau,r_0) d\tau \right). $$
This satisfies 
$$ \begin{array}{rl}\lambda^{-1}d\lambda= &  -ih(x,y,t,s) ds
-if(x,y,t,r_0) dt \\[2mm]
& + (-i \int_{r_0}^s \partial_t
h(x,y,t,\sigma) d\sigma) dt \\[2mm]
&  -i d_{T^2} \left( \int_{r_0}^s
h(x,y,t,\sigma)d\sigma + \int_0^t f(x,y,\tau,r_0) d\tau \right) \\[2mm]
= & -ih(x,y,t,s) ds -i f(x,y,t,s)dt  \\[2mm] & -i d_{T^2}\left(
\int_{r_0}^s h(x,y,t,\sigma)d\sigma + \int_0^t f(x,y,\tau,r_0) d\tau
\right), \end{array} $$
where the equality follows from the last equation of (\ref{gaugeF}).
The gauge transformed connection $\lambda \A$ is of the form
$$ \lambda \A= a_{T^2}-i d_{T^2} \left( \int_{r_0}^s
h(x,y,t,\sigma)d\sigma + \int_0^t f(x,y,\tau,r_0) d\tau \right)  =
\tilde a_{T^2}. $$  
This satisfies the equations
$$ \begin{array}{c} F_{\tilde a_{T^2}}=0 \\
\partial_t \tilde a_{T^2}=\partial_t a_{T^2} -id_{T^2}(f(x,y,t,r_0) +
\int_{r_0}^s \partial_t h(x,y,t,\sigma) d\sigma) = \\
=\partial_t a_{T^2} -i d_{T^2}f =0 \\
\partial_s \tilde a_{T^2}=\partial_s a_{T^2} -i d_{T^2}h =0. \end{array} $$
This completes the proof of the Claim. 

The asymptotic limits $a_\infty^\pm$ of
$(A^\prime(\pm\infty),\psi^\prime(\pm\infty))$ then satisfy
\[ [ a_\infty^+ ]=[ a_\infty^- ]=[a_{T^2}] \]
in the universal cover $\chi_0(T^2,Y)$ of $\chi(T^2)$. 
This completes the proof of the Proposition.

\endproof

Notice, however, that in general the uniform convergence of
Proposition \ref{unif:conv} should not be expected. In fact, we have
seen that the
estimate (\ref{exp:decay}) and Lemma \ref{decay:bound} give
$$ \| F_{A_r(t)} \|^2_{L^2(T^2\times [r_0,r_1]\times \{ t \})} \leq C
e^{-\frac{c|t|}{r}}. $$
This allows for non-uniformity of the convergence: there may be 
sequences $T_r \to \infty$ for which (\ref{uniform:curv}) is violated,
and an estimate like 
$$ \| F_{A_r(t)} \|^2_{L^2(T^2\times [r_0,r_1]\times \{ t \})}
\geq c >0 $$ 
holds
for all $|t|\in [T_r -\ell/2, T_r +\ell/2]$, for some $\ell >0$ 
independent of $r\geq r_0$. This is clearly compatible with the
constraint (\ref{exp:decay}).

In the next subsection we analyze the geometric limits of monopoles
$(\A_r,\Psi_r)$ on the domains $T^2\times [-r,r]\times [T_r,\infty)$
and $T^2\times [-r,r]\times (-\infty,-T_r]$. This is where we
encounter the non-uniform convergence which is not detected by limits
on compact sets.

\subsection{Holomorphic disks in $\chi(T^2)$}

Consider the Seiberg-Witten equations in
on $T^2\times [-r,r]\times \R$. Let $(\A_r,\Psi_r)$
be a family of solutions, where we write the connections as
$$ \A_r=a_r(w,s,t)+f_r(w,s,t)ds +h_r(w,s,t)dt $$ 
and the spinors as $\Psi_r= (\alpha_r(w,s,t),\beta_r(w,s,t))$. 

Recall again that gauge equivalent solutions $(A_r(t),\psi_r(t))$, in
a temporal gauge in the $t\in \R$ direction, decay exponentially at
the asymptotic ends $t\to\pm\infty$ to critical points
$(A_r(\pm\infty),\psi_r(\pm\infty))$ of the Chern-Simons-Dirac
functional. 

Suppose given a fixed compact set $T^2\times [r_0,r_1]$ in $Y(r)$. 
Choose a sequence $T_r\to\infty$, for which the estimate
(\ref{exp:decay}) gives
\ba \| a_r(w,s,\pm T_r)- a^\pm (w) \|^2 + \| f_r(w,s,\pm
T_r)-f^\pm(w,s) \|^2 \leq C e^{-c r}, \label{asympt:close} \na
in the $L^2_1$-norm on $T^2\times [r_0,r_1] $,
where $a^\pm$ are the limits of $(A_r(\pm\infty),\psi_r(\pm\infty))$
as $r\to\infty$ on $T^2\times [r_0,r_1] $. 
In the first part \cite{CMW} of this work we showed that, in the
gluing of solutions to the 3-dimensional Seiberg-Witten equations
(\ref{3SW}), we are only 
interested in critical points $(A_r(\pm\infty),\psi_r(\pm\infty))$ of
the Chern-Simons-Dirac functional that break through smooth points
$a^\pm\neq \vartheta$ in the character variety $\chi(T^2)$, as
$r\to\infty$. Thus, we can 
assume that the elements $f^\pm(w,s)$
are exponentially small on $T^2\times [r_0,r_1] $, for large enough
$r_0>0$, 
$$ \| f^\pm(w,s) \|^2_{L^2(T^2\times \{ s \})} \leq  C e^{-\delta s}. $$
This exponential estimate follows from the analysis of the asymptotics 
of monopoles on a 3-manifold with an end modeled on $T^2\times
[0,\infty)$, as in \cite{CMW}. 

This simple observation implies that we should expect more complicated
geometric limits than the limits on compact sets analyzed in the
previous section. In fact, for $T_r$ growing sufficiently fast, we have
convergence to the limits $a^\pm$ on
$T^2\times [r_0,r_1]\times \{ t \}$, with $|t| \geq T_r$. Thus, whenever
$a^+ \neq a^-$, we must
have geometric limits which are not just constant flat connections on
$T^2$. 

As discussed
previously, in polar coordinates 
$z=s+it$ and $z=e^{\rho +i\theta}$, with
$$ a_r(w,\rho,\theta)=a_r(w,e^{\rho +i\theta}) $$
$$ f_r(w,\rho,\theta)=e^{-\rho} \cos\theta \ h_r(w,e^{\rho
+i\theta})-e^{-\rho}\sin\theta \ f_r(w,e^{\rho +i\theta}) $$
$$ h_r(w,\rho,\theta)=e^{-\rho}\cos\theta \ f_r(w,e^{\rho+i\theta})
+e^{-\rho}\sin\theta \ h_r(w,e^{\rho+i\theta}), $$
the curvature equation takes the form
$$ \partial_\rho a_r -dh_r =*(\partial_\theta a_r -df_r +
i(\bar\alpha_r \beta_r +\alpha_r\bar\beta_r)) $$
$$ \partial_\rho f_r -\partial_\theta h_r=e^{2\rho} *( F_{a_r} +
\frac{i}{2}(|\alpha_r|^2-|\beta_r|^2)\omega). $$

We are now interested in studying the convergence away from
arbitrarily large compact sets. That is, we are interested in the
convergence on regions of the form
$$ \begin{array}{rl} \tilde\Omega(r)= & 
T^2\times [0,r] \times [T_r,\infty) \\ & \cup T^2\times [0,r) \times
(-\infty, -T_r] \\ & \cup T^2 \times \{ \rho \in [\log T_r,\infty) \
\theta\in [-\pi/2,\pi/2] \}\end{array}. $$

We choose a sequence $T(r)$ such that, as in (\ref{asympt:close}), the
finite energy solutions $(\A_r,\Psi_r)$ on $Y(r)\times \R$ are
exponentially close to the asymptotic values $a^\pm$ on 
$$ T^2\times [-r_0,r_0]\times ([T(r),\infty)\cup (-\infty,-T(r)]. $$
We write $T(r)=e^{T_r}$, with $T_r\to\infty$ as $r\to\infty$.

In order to study the convergence on the cylinder $T^2\times [-r,r]\times
\R$, we introduce a suitable
rescaling of the coordinates on $T^2\times \R^2$.
In polar coordinates, we introduce the translation $\rho \mapsto \rho-
T_r$. In these new 
coordinates, the second part of the curvature equation becomes
$$ \partial_\rho f_r -\partial_\theta h_r -e^{2T_r} e^{2\rho} *( F_{a_r} +
\frac{i}{2}(|\alpha_r|^2-|\beta_r|^2)\omega)=0. $$

Consider the unit disk $D=\{ \rho\leq 0 \}$ in the new coordinates. By
the previous analysis 
of the convergence on compact sets we obtain smooth convergence of a
subsequence to a solution of
$$ \partial_\rho a -*(\partial_\theta a-df) =0 $$
$$ F_a =0. $$
Up to a gauge transformation, this is a holomorphic map
$$ a: D \to \chi_0(T^2,Y)=\C. $$
Namely, $a(\rho,\theta)$ is the flat connection on $T^2$ satisfying
$$ a(\rho,\theta)(w)=a(w,\rho,\theta). $$
With a slight abuse of notation we identify $a(\rho,\theta)$ with the
class in $\chi_0(T^2,Y)$ obtained modulo homotopically trivial gauge
transformations.  
The map $a$ has
the property that it maps the center of the disk, $\rho\to
-\infty$, to the flat connection $a_\infty$, which is the limit on
compact sets of the original sequence $(\A_r,\Psi_r)$. Moreover, by
the finite energy condition, the map $a$ is holomorphic up to the
boundary, and the restriction to the boundary
$$ a: S^1=\partial D = \{ \rho=0 \} \to \chi_0(T^2,Y)=\C $$
gives  a curve in $\chi_0(T^2,Y)$. By our choice of $T_r$, this curve
will satisfy $a(0,\pi/2)=a^+$ and $a(0,-\pi/2)=a^-$.

Thus, we have obtained that the convergence on the $T^2\times
[-r,r]\times \R$ part of $Y(r)\times \R$, away from compact sets,
produces a holomorphic map. Now we want to describe the corresponding
non-uniform convergence, away from compact sets, on the 
sub-domains of $Y(r)\times \R$ that correspond to the knot complement
and to the tubular neighborhood of the knot. We will use a
similar rescaling of the coordinates.

Let $V_r$ be the knot complement $V$ with a cylinder of length $r$,
$$ V_r = V\cup_{T^2\times \{ 0 \}} T^2\times [0,r]. $$

On the manifold $V_r\times \R$ consider the domain 
$$ \tilde\Omega(r) = V_r\times ([T_r,\infty) \cup (-\infty, -T_r]) $$
$$ \cup T^2\times \{ \rho\in [\log T_r, \infty) \
\theta\in[-\pi/2,\pi/2] \}, $$  
where we use the polar coordinates
$(\rho,\theta)$ introduced before.
On this domain consider the change of coordinates
$$ \rho\mapsto \rho - T_r $$
in the $\rho\in \R$ direction, and the 
corresponding change of coordinates
$$ |t|\mapsto e^{-T_r} |t| $$
in the $t\in \R$ direction, where the coordinates $\rho$ and $t$ are
related by  $z=s+it$ and
$z=e^{\rho +i\theta}$, as before, with $(s,t)\in [0,\infty)\times \R$.

This produces the rescaled domain
$$ V_r\times ([T_r e^{-T_r},\infty) \cup (-\infty, -T_r e^{-T_r}]) $$
$$ \cup T^2\times \{ \rho\in [\log(T_r) - T_r , \infty) \
\theta\in[-\pi/2,\pi/2] \}. $$ 

\begin{figure}[ht]
\epsfig{file= 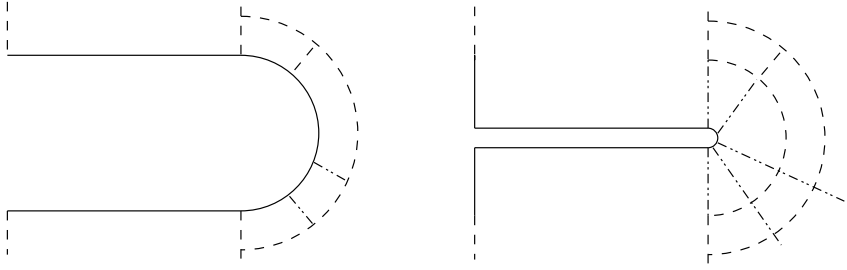}
\caption{The domain $\tilde\Omega(r)$ and its rescaling
\label{figII1}} 
\end{figure}

Consider the domains 
$$ \Omega(r)= V_r\times ([-1,-T_r e^{-T_r}]\cup[T_r e^{-T_r},1]) $$
$$ \cup   T^2\times \{ \rho\in 
(-T_r +\log(T_r),0] \ \theta\in[-\pi/2,\pi/2] \}, $$
and
$$ \Omega_1(r)= V_r\times ([-1,-T_r e^{-T_r}]\cup[T_r e^{-T_r},1])
\subset \Omega(r), $$
$$ \Omega_2(r)= T^2\times \{ \rho\in 
(-T_r +\log(T_r),0] \ \theta\in[-\pi/2,\pi/2] \}
\subset \Omega(r). $$ 
The finite energy solutions $(\A_r,\Psi_r)$ of the SW equations
satisfy
$$ e^{-T_r} \frac{\partial}{\partial_t}A_r(t)= -*F_{A_r(t)}+
\sigma(\psi_r(t),\psi_r(t)) $$
$$ e^{-T_r} \frac{\partial}{\partial_t} \psi_r(t)=
-\partial_{A_r(t)}\psi_r(t) $$ 
on $\Omega_1(r)$. 

\begin{figure}[ht]
\epsfig{file= 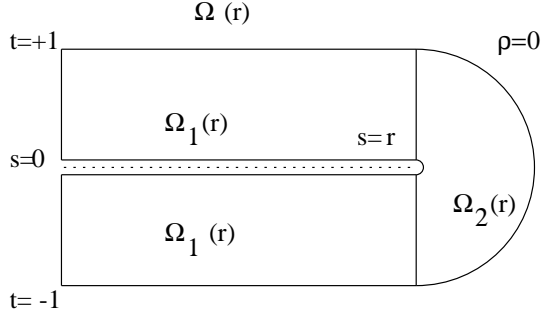,angle= 270}
\caption{The domains $\Omega(r)$, $\Omega_1(r)$, and $\Omega_2$
\label{figII2}} 
\end{figure}

The finite energy condition ensures that the left hand side converges
to zero uniformly on any fixed compact set. Thus, the family
$(\A_r,\Psi_r)$, under the present change of coordinates, converges
smoothly on compact sets in to a path $(A(t),\psi(t))$ of finite
energy solutions in the domain
$V \times \R^*$, where $\R^*=\R \backslash \{ 0 \}$, where   
$V$ is endowed with an infinite cylindrical end $T^2\times
[0,\infty)$. That is, we obtain a path $(A(t),\psi(t))$ of elements 
satisfying 
$$ {}*F_{A(t)}+
\sigma(\psi(t),\psi(t))=0 $$
$$ \partial_{A(t)}\psi(t)=0. $$
This defines a path $[A(t),\psi(t)]$ in ${\cal M}_V$ for $0< t \leq
1$. The endpoint satisfies
$$ [A(1),\psi(1)]=[A_{+\infty}',\psi_{+\infty}'], $$ 
where 
$(A_{+\infty}',\psi_{+\infty}')$ is the limit on $V$ as $r\to\infty$ of
the asymptotic values $(A_r(+\infty),\psi_r(+\infty))$ of
$(\A_r,\Psi_r)$. Thus, we have 
$$ \partial_\infty(A_{+\infty}',\psi_{+\infty}')=a^+, $$
under the boundary value map.
Moreover, we obtain an element 
$$ \lim_{t\to 0_+}[A(t),\psi(t)]=[A,\psi]$$ 
in ${\cal M}_V$, satisfying
$$ \partial_\infty(A,\psi)= a_\infty', $$
where $a_\infty'$ is the asymptotic value of the limit on compact sets
of the solutions $(\A_r,\Psi_r)$ on $V\times \R$, as analysed in the
previous section. 
Similarly, we obtain a path $[A(t),\psi(t)]$ in ${\cal M}_V$ for
$-1\leq t < 0$, where the endpoint satisfies
$$ [A(-1),\psi(-1)]=[A_{-\infty}',\psi_{-\infty}'], $$ 
where 
$(A_{-\infty}',\psi_{-\infty}')$ is the limit on $V$ as $r\to\infty$ of
the asymptotic values $(A_r(-\infty),\psi_r(-\infty))$ of
$(\A_r,\Psi_r)$. Thus, we have 
$$ \partial_\infty(A_{-\infty}',\psi_{-\infty}')=a^-, $$
under the boundary value map.
Similarly, the element 
$$ \lim_{t\to 0_-} [A(t),\psi(t)]=[\tilde A,\tilde \psi]$$ 
in ${\cal M}_V$ satisfies 
$$ \partial_\infty(\tilde A,\tilde\psi)= a_\infty'. $$
 
Notice that in general the element $[A,\psi]$ and the element $[\tilde
A,\tilde \psi]$ need not be the same element in ${\cal M}_V$. In fact,
the domain $\Omega_1(r)$ is defined by rescaling the initial domain
$\tilde\Omega(r)$, and the values of solutions near the frontier
$t=\pm T_r$ of $\tilde\Omega(r)$ need not approach the same limits as
$r\to\infty$. 

{}From the information on the convergence on compact sets
obtained in the previous section, together with the choice of the
domains, we can only derive the relation
$$ \partial_\infty(A,\psi)=\partial_\infty(\tilde A,\tilde\psi). $$
Thus, both $[A,\psi]$ and $[\tilde A,\tilde \psi]$ are in the
zero-dimensional fiber $\partial_\infty^{-1}(a_\infty ')$,
cf. \cite{CMW}. 

Moreover, by the previous analysis, together with the analysis of
asymptotics on $V\times \R$ of the previous section, we obtain
convergence on the $\Omega_2(r)$ to a function $a: D^+\to
\chi_0(T^2,V)$, holomorphic in a neighborhood $U$ of the half disk 
$$ D^+=\{ \rho\in (-\infty,0]\ \theta\in [-\pi/2,\pi/2] \}, $$ 
which maps the center of 
the disk to $a_\infty'$, and such that the restriction to the boundary
component $\{ \theta \in \{ -\pi/2,\pi/2 \} \}\subset \partial D^+$ gives a
curve in $\partial_\infty({\cal M}_V) \subset \chi_0(T^2,V)$ connecting
$a^+$, $a_\infty'$,  and $a^-$. The image of the component $\{ -\pi/2
\leq \theta \leq \pi/2 , \ \rho=0 \}$ of $\partial D^+$ is mapped to
another curve in $\chi_0(T^2,V)$, not necessarily contained in
$\partial_\infty({\cal M}_V)$, which connects the points $a^+$,
$a_\infty'$,  and $a^-$. The image of $\partial D^+$ is a closed curve
in $\chi_0(T^2,V)$, not necessarily smooth across the points
$\theta=\pm\pi/2$.

\subsection{The convergence theorem}

We need one more preliminary Lemma, which connects the limits on
compact sets, with the non-uniform limits obtained in the previous
section.   

\begin{Lem}
Let $(\A,\Psi)$ be a finite energy solution of the Seiberg--Witten
equations (\ref{4SW}) on the 4-manifold $V\times \R$ with the infinite
end $T^2\times [0,\infty)\times \R$. Let $a_\infty \in \chi_0(T^2,V)$
be the radial gauge limit of $(\A,\Psi)$. Let $(A(t),\psi(t))$ be
gauge equivalent to the original $(\A,\Psi)$, in a temporal gauge in
the $t\in\R$ direction. Then the elements $(A(t),\psi(t))$ converge
in the $s\in [r_0,\infty)$ direction to $a_\infty$, uniformly with
respect to $|t|\geq T_0$, and we have limits
$$ \lim_{t\to +\infty} (A(t),\psi(t))=(A,\psi), $$
and 
$$ \lim_{t\to -\infty} (A(t),\psi(t))=(\tilde A,\tilde\psi), $$
with
$$ \partial_\infty(A,\psi)=\partial_\infty(\tilde
A,\tilde\psi)=a_\infty. $$
\label{ts:limits}
\end{Lem}

\proof
In order to prove the convergence in the time direction, consider the
asymptotics of solutions in radial gauge studied in the previous
section. From the analysis of the system of ODE's
(\ref{SWradfourier}), we obtained rates of decay as $\rho\to
\infty$ of the form $\exp(-n \rho)$ or $\exp(n\rho - C_{nlk} e^\rho)$,
or $\exp(- C_{nlk} e^\rho)$, for some positive constants $C_{nlk}>0$. 

In the last case, at a
point $(s,t)$ along the end $[r_0,\infty)\times \R$, the $L^2$-norm on
$T^2$ is bounded by $\exp(-C (s^2+t^2)^{1/2})$. In particular,
we have exponential decay in the $s\in
[r_0,\infty)$ direction, uniformly in $|t|\geq T_0$, bounded by
$\exp(-C T_0 s)$. Thus, if we consider a temporal gauge
representative $(A'(t),\psi'(t))$ of $(\A',\Psi')$, for every sequence
$|t_k| \to\infty$ we obtain convergence up to gauge of a subsequence
to an element in $\M_V$, exponentially decaying in the $s\in
[r_0,\infty)$ direction to the asymptotic value $a_\infty'$, at a rate
at least equal to $\exp(-|\lambda| T_0 s)$. 

The case of the slow decaying solutions with asymptotics $e^{-n \rho}$
gives polynomial decay in the $s\in [r_0,\infty)$ direction for every
fixed $t\in\R$, at a rate 
$$ \frac{1}{\left(1+\left(\frac{t}{s}\right)^2 \right)^{n/2}} \cdot
\frac{1}{s^n} \leq  \frac{1}{s^n}, $$
uniform in $|t|\geq T_0$. Thus, in this case, the temporal gauge
representative $(A'(t),\psi'(t))$ converges to solutions of
(\ref{3SW}) on $V$, as $t\to\pm\infty$, that decay to the asymptotic
value  $a_\infty'$ at a polynomial rate $1/s^n$. The analysis of
\cite{CMW} on the center manifold for the equations (\ref{3SW}) on
$T^2 \times [r_0,\infty)$ (cf. \cite{MMR}), imply the following. If
the decay is asymptotic to $1/s^n$, with $n\geq 2$, then the actual
rate of decay of the limit solution is exponential and we have
$$ a_\infty' \in \partial_\infty(\M_V^*). $$
If the rate of decay is $1/s$, then the limit solution might also
decay like $1/s$ along $s\in [r_0,\infty)$. In this case the
asymptotic value is $ a_\infty' =\vartheta$, the unique ``bad point''
in the character variety $\chi(T^2)$.

\endproof

We can summarize all the previous results as follows

\begin{The}
On the manifold $Y(r)\times \R$ consider classes of finite energy solutions
$[\A_r,\Psi_r]$ of (\ref{4SW}), with asymptotic values 
$$ [A_r(\pm\infty),\psi_r(\pm\infty)]=
[A_{\pm\infty}',\psi_{\pm\infty}'] \#^r_{a^\pm} 
[a^\pm,0], $$
as $t\to \pm\infty$. Under the usual assumptions on the
scalar curvature on $\nu(K)$, we obtain the following geometric limits
as $r\to\infty$.

(i) a finite energy $[\A',\Psi']$ on $V\times \R$, with radial 
limit $a_\infty' \in \chi_0(T^2,V)$, and limits as $t\to\pm\infty$
$$ [A,\psi] \ \hbox{and} \ [\tilde A, \tilde \psi] \ \hbox{in} \
\partial_\infty^{-1}(a_\infty') \subset \M_V. $$

(ii) a solution $[a_\infty '',0]$ on $\nu(K)\times \R$, with $a_\infty
''\in \chi_0(T^2,\nu(K))$ .

(iii) a path $[A(t),\psi(t)]\in \M_V$  for $t\in [-1,0)\cup (0,1]$
with 
$$ [A(\pm 1),\psi(\pm 1)]=[A_{\pm\infty}',\psi_{\pm\infty}'] $$
and the limits as $t\to 0_\pm$ equal to $[A,\psi]$ and $[\tilde A,
\tilde \psi]$.

(iv) a path $[A(t),0]$ in $\M_{\nu(K)}=\chi(\nu(K))$ for $t\in [-1,1]$
satisfying $[A(\pm 1),0]=[a^\pm, 0]$ and $[A(0),0]=[a_\infty '',0]$. 

(v) a holomorphic map $a: D\to \chi_0(T^2,Y)$, with $a(0)=a_\infty$,
the limit on compact sets on the cylinder $T^2\times [-r,r]\times
\R$. 

(vi) maps $a_V: D^+ \to \chi_0(T^2,V)$ and $a_\nu: D^+\to
\chi_0(T^2,\nu(K))$, holomorphic on a neighborhood of the half disk
$D^+$ with $a_V(\rho,\pm\pi/2)=[A(t),\psi(t)]$ and
$a_\nu(\rho,\pm\pi/2)=[A(t),0]$. 

\label{limits}
\end{The}

\proof
Suppose given a family $(\A_r,\Psi_r)$ of finite energy solutions of
the Seiberg-Witten equations (\ref{4SW}) on the manifolds $Y(r)\times
\R$, with asymptotic values $(A_r(\pm\infty),\psi_r(\pm\infty))$ as
$t\to\pm\infty$ satisfying (\ref{3SW}). Assume that, as $r\to\infty$
the asymptotic values split as
$$ (A_r(\pm\infty),\psi_r(\pm\infty))=
(A_{\pm\infty}',\psi_{\pm\infty}') \#^{r}_{a^\pm}\lambda
(a^\pm,0), $$ for some gauge transformation $\lambda$ on $\nu(K)$,
according to \cite{CMW}.

Then, up to gauge, there is a subsequence, which we still denote
$(\A_r,\Psi_r)$, with the following behavior.
On compact sets contained in the long cylinder $T^2\times [-r,r]\times 
\R$ it converges smoothly to a flat
connection $a_\infty$ on $T^2\times \R^2$, constant in the $\R^2$
directions.
On compact sets in $V\times \R$ it converges smoothly to a
finite energy solution  
$(\A',\Psi')$ of the Seiberg-Witten equations (\ref{4SW}) on $V\times \R$,
which is exponentially decaying in radial gauge on the end $T^2\times
[0,\infty)\times \R$ to a flat connection $a_\infty'$.
This follows from the analysis of the asymptotics in the previous
section. 

The solution $(\A',\Psi')$ in a temporal gauge in the $t\in\R$ direction 
converges to solutions $(A,\psi)$ and $(\tilde A,\tilde \psi)$ of the
Seiberg-Witten equations (\ref{3SW}) on $V$, as $t\to\pm\infty$
respectively. The elements $(A,\psi)$ and $(\tilde A,\tilde \psi)$
have asymptotic value along the end $T^2\times \R$ of $V$
$$ \partial_\infty(A,\psi)=\partial_\infty(\tilde A,\tilde
\psi)=a_\infty', $$
under the map $\partial_\infty$ analyzed in \cite{CMW}. This result
follows from Lemma \ref{ts:limits}

Finally, again from the arguments presented in the previous sections,
on compact sets in $\nu(K)\times \R$ the subsequence $(\A_r,\Psi_r)$
converges to a solution of the 
abelian ASD equation, that is, up to gauge, to a flat connection 
$a_\infty ''$ on $T^2$.

Up to the changes of coordinates $\rho\mapsto \rho- T_r$
and $|t|\mapsto e^{-T_r} |t|$, for a sufficiently fast growing
$T_r\to\infty$, there is a subsequence of
$(\A_r,\Psi_r)$ that converges smoothly on the domain $\Omega_1$
defined above to a path $[A(t),\psi(t)]$ in ${\cal M}_V$ for $t\in
[-1,0)\cup (0,1]$ with 
$$[A(-1),\psi(-1)]=[A_{-\infty}',\psi_{-\infty}'], \ \ \hbox{with} \ \ 
\partial_\infty(A_{-\infty}',\psi_{-\infty}')= a^- $$ 
$$ \lim_{t\to 0_-} [A(t),\psi(t)]=[\tilde A,\tilde \psi],\ \ \hbox{with} \ \ 
\partial_\infty(A,\psi)=a_\infty ' $$
$$ \lim_{t\to 0_+} [A(t),\psi(t)]=[A,\psi], \ \ \hbox{with} \ \ 
\partial_\infty(A,\psi)=a_\infty ' $$
$$[A(1),\psi(1)]=[A_{+\infty}',\psi_{+\infty}']\ \ \hbox{with} \ \ 
\partial_\infty(A_{+\infty}',\psi_{+\infty}')= a^+.$$ 
Under the same change of coordinates there is a subsequence
that converges smoothly on the domain $\Omega_2$ described above to a
function $a_V: D^+\to \chi_0(T^2,Y)=\C$, holomorphic in a
neighborhood $U$ of the half disk $D^+$,
which maps the center of 
the disk to $a_\infty'$ and such that the restriction to the boundary
component $\theta \in \{ -\pi/2, \pi/2 \}$
of $D^+$ gives a curve in $\chi_0(T^2,Y)$
connecting $a^+$, $a_\infty'$, and $a^-$. 

Similarly, on the manifold $\nu(K)\times \R$ we can consider similarly 
defined domains $\Omega_1$ and $\Omega_2$, and an analogous
reparametrization. We obtain 
convergence on $\Omega_1$ to a path $[A(t),0]$ in ${\cal
M}_{\nu(K)}=\chi(\nu(K))$ for $t\in [-1,1]$ with 
$$[A(-1),0]=[a^-,0], $$
$$[A(0),0]=[a_\infty '',0], $$
$$[A(1),0]=[a^+,0], $$
where $a_\infty ''$ is the limit on compact sets of the
$(\A_r,\Psi_r)$ on $\nu(K)\times \R$, up to gauge.
On $\Omega_2$ we obtain convergence to a function $a_\nu:
D^+\to \chi_0(T^2,Y)$, holomorphic in a
neighborhood $U$ of the half disk $D^+$, which maps the center of 
the disk to $a_\infty ''$, and such that the restriction to the
boundary component $\theta \in \{ -\pi/2, \pi/2 \}$
of $D^+$ gives a curve in $\chi_0(T^2,Y)$
connecting $a^+$, $a_\infty ''$, and $a^-$.

Finally, the reparametrization $\rho\mapsto \rho- T_r$ on the long
cylinder $(s,t)\in [-r,r]\times \R$ gives convergence on the disk
$D=\{ \rho\leq 0 \}$ to a holomorphic map 
$ a: D\to \chi_0(T^2,Y)$, which maps the center of 
the disk to $a_\infty$. 

There are subdomains $D_1^+$ and $D_2^+$ inside the disk $D$, and
conformal equivalences $\varphi_i : D^+ \to D_i^+$, such that we have
$$ a|_{D_1^+} = a_V \circ \varphi_1^{-1} $$
$$ a|_{D_2^+} = a_\nu \circ \varphi_2^{-1}. $$

In the special case where $a^+=a^-=a_\infty$ up to gauge, the holomorphic
functions are constant and the paths $[A(t),\psi(t)]$ in ${\cal M}_V$
and ${\cal M}_{\nu(K)}$ are also constant.

\endproof

Notice that, in general, the points $a_\infty'$, $a_\infty ''$, and
$a_\infty$ are distinct points in $\chi_0(T^2,Y)$. The point
$a_\infty'$ is constrained to be on the 1-dimensional subspace
$\partial_\infty(\M_V)$ and the point $a_\infty ''$ is contained in
$\M^{red}_{\nu(K),\mu}= \chi(\nu(K))$. 

We have seen that, if the three
points $a^+$, $a^-$, and $a_\infty$ coincide, $a^+=a^-=a_\infty$ up to
gauge, then the flow line 
$(\A_r,\Psi_r)$ on $Y(r)\times \R$ must connect critical points which
break through the same asymptotic values, $a^+=a^-=a_\infty$.

We can also observe that the following Corollary holds.

\begin{Cor}
Assume that the points $a_\infty'$ and $a_\infty''$ coincide up to
gauge. Assume, moreover, that they are distinct from the bad point
$\vartheta$, and from both $a^+$ and $a^-$, 
$$ a^+\neq a_\infty' \neq a^-. $$ then the original
flow line $(\A_r,\Psi_r)$ on $Y(r)\times \R$ must connect critical
points with relative index at least two.
\end{Cor}

\proof Let $(\A',\Psi')$ be the limit on compact sets in $V\times \R$
of the family $(\A_r,\Psi_r)$ as $r\to\infty$. Then the point
$a_\infty'$ is the radial limit of $(\A',\Psi')$. By Lemma
\ref{ts:limits} we know that as $t\to \pm\infty$ the solution
$(\A',\Psi')$ has asymptotic values $[A, \psi]$ and $[\tilde A, \tilde
\psi]$ in $\M_V^*$ with
$$ \partial_\infty[A, \psi]=\partial_\infty[\tilde A, \tilde
\psi]=a_\infty. $$
By gluing this solution with the solution $(a_\infty'',0)$ on the
$\nu(K)\times \R$ side, we obtain a flow line in $Y(r)\times \R$ which
connects the critical points 
$$ [A, \psi]\#^r_{a_\infty} [a_\infty, 0] $$
and
$$ [\tilde A, \tilde\psi]\#^r_{a_\infty} [a_\infty, 0]. $$
This shows that there are two intermediate critical points between the
endpoints of the flow line $(\A_r,\Psi_r)$.

\endproof

\section{Gluing}

In this section we describe how to produce an approximate solution to
the 4-dimensional Seiberg--Witten equations on $Y(r(T))\times \R$,
with gluing parameter $T\geq T_0$, by pasting together the various
geometric limits of Theorem \ref{limits}.

\subsection{The approximate solutions}

We define the following quantities that will be used throughout this
section:
$$ R(T)=\frac{1}{\pi}( e^T + T ), $$
$$ \ell(T)= \frac{1}{\pi}( e^T + T ) \sin \left( \frac{\pi T}{e^T + T}
\right), $$
$$ r(T)=  \left( r_0 + \frac{1}{\pi}( e^T + T ) \cos\left( \frac{\pi
T}{e^T + T} \right)\right). $$
Moreover, we define
$$ \Upsilon (t,T)= r(T)-( R(T)^2 -t^2 )^{1/2}, $$
for $|t|\leq R(T)$, and
$$ \tilde \Upsilon(t,T)= r_0 + \cos\left( \frac{\pi
T}{e^T + T} \right) - \min\left\{ (1-t^2)^{1/2}, \cos\left( \frac{\pi
T}{e^T + T} \right) \right\}, $$
for $|t|\leq 1$.

Consider the path $[A(t),\psi(t)]$ on the domain
\ba V_{r(T)}\times ([-1,-\ell(T)/R(T)]\cup [\ell(T)/R(T),1],
\label{dom0} \na  
with $T$ large enough, so that the elements $[A(\pm \ell(T)/R(T)),\psi(\pm
\ell(T)/R(T))]$ are very close to the elements $[A,\psi]$ and $[\tilde A,
\tilde \psi]$, and with $r_0$ large enough, so that at $s=r_0$ the
elements $[A(t),\psi(t)]$ are very close 
to the asymptotic values $a_\infty(t)$, up to an error of the order of
$e^{-\delta r_0}$.
Consider the rescaled domain
\ba V_{r(T)}\times ([-R(T),-\ell(T)]\cup [\ell(T),R(T)]), \label{dom1}
\na 
under the change of coordinates $|t|\to R(T) |t|$,
and write $[A(t),\psi(t)]_T$ for the corresponding rescaled element.
Consider then the path $[A(t),\psi(t)]_T$ restricted over the domain
\ba  {\cal R}_1(T)= \begin{array}{l} (V\times I(T))\cup_{T^2\times \{
s=0 \}\times 
I(T)} \\
 T^2\times \{ (s,t) \ | \  s\in [0,\Upsilon(t,T)], \ \ t\in
I(T)  \}, \end{array} \label{R1} \na 
with
\ba \label{intervalT} I(T)=[-R(T),-\ell(T)]\cup [\ell(T),R(T)], \na
inside the rescaled domain (\ref{dom1}).

By Theorem \ref{limits}, we can consider a disk $D_T$ of radius $R(T)$,
and a map $a_T: D_T\to \chi_0(T^2,Y)$, such that the center of
the circle is mapped to the limit $a_\infty$ of Theorem \ref{limits},
and the points along the boundary 
corresponding to the angles $\{ -\pi, -\pi/2, 0, \pi/2 \}$ are mapped
to the points $\{ a_\infty '', a^-, a_\infty ', a^+ \}$,
respectively. The function $a_T: D_T\to \chi_0(T^2,Y)$ agrees with the
map $a: D\to \chi_0(T^2,Y)$ of Theorem \ref{limits} restricted to a
subdomain of the unit disk $D$ homeomorphic to $D_T$. We can
guarantee, by choosing $T$ large enough, that on the arc of length
$2T$ on $\partial D_T$, centered at the angle $\theta=0$, the values
of the function $a_T$ are sufficiently close to $a_\infty '$, and on
the similar arc centered at the angle $\theta=-\pi$ the values
of the function $a_T$ are sufficiently close to $a_\infty ''$.

Consider then the element $(\A',\Psi')$ restricted over the domain 
\ba {\cal R}_3(T)=V_{r_0} \times [-\ell(T),\ell(T)], \label{R3} \na
For $T\geq T_0$ large enough, we have that, near $t=\pm \ell(T)$,
$(\A',\Psi')$ is sufficiently close to the elements $[A,\psi]$ and
$[\tilde A, \tilde \psi]$, and, near $s=r_0$, we have that
$(\A',\Psi')$ is sufficiently close to the asymptotic value $a_\infty$.

Moreover, consider the rescaled path $(a ''(t), 0)_T$ on the domain
\ba \label{R4}  {\cal R}_4(T)=\begin{array}{l} (\nu(K)\times
I(T))\cup_{T^2\times \{ s=0 
\}\times I(T)}  \\ 
 T^2\times \{ (s,t) \ | \  s\in [0,\Upsilon(t,T)], \ \ t\in
I(T) \} , \end{array} \na
with $I(T)$ as in (\ref{intervalT}).

Consider then the element $(a_\infty '', 0)$ on the domain
\ba \label{R5} {\cal R}_5(T)=\nu(K)_{r_0}\times
[-\ell(T),\ell(T)],  \na 

\begin{figure}[ht]
\epsfig{file= 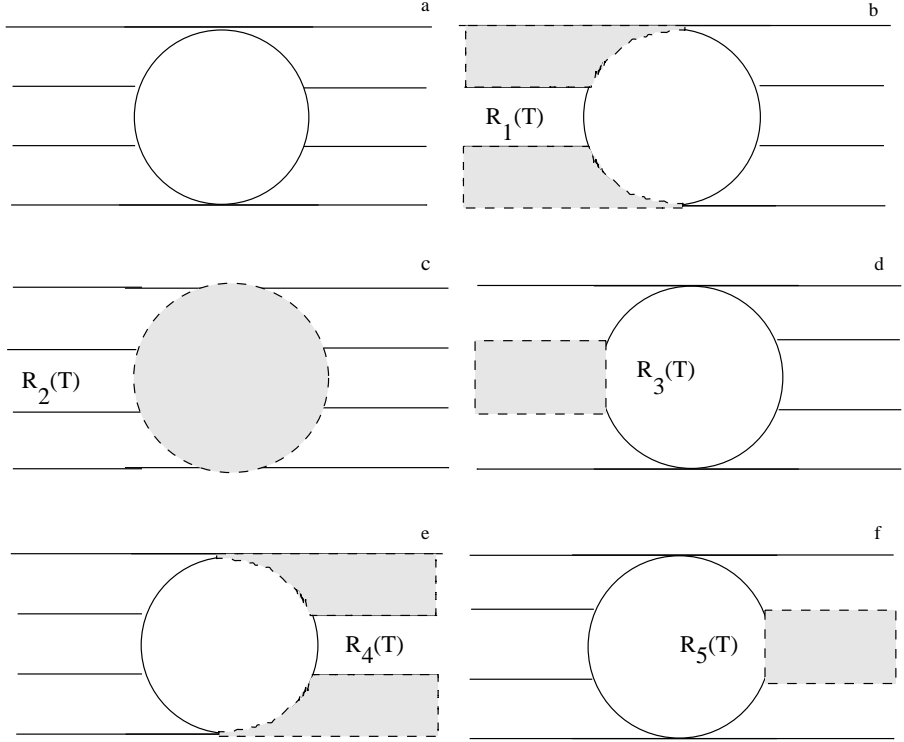}
\caption{The construction of the approximate solution
\label{figII4}} 
\end{figure}

Now consider the manifold $Y(r(T))\times \R$, with the long cylinder
$T^2\times [-r(T),r(T)]\times \R$.
Consider the product region 
\ba T^2\times \{ s\in [-r(T),r(T)], \ \  t\in [-R(T),R(T)] \}
\label{rectangle} \na  in $Y(r(T))\times \R$
as illustrated in part {\em (a)} of Figure \ref{figII4}.

Consider inside (\ref{rectangle}) the region ${\cal R}_2(T)$ given by
$T^2\times D_T$, with the
disk $D_T$ centered at $s=0$, $t=0$, as in Figure
\ref{figII4}. We identify the regions ${\cal 
R}_i(T)$ of (\ref{R1}), (\ref{R3}), (\ref{R4}), (\ref{R5}), with the
corresponding regions as in Figure  
\ref{figII4}. The smaller strip of regions ${\cal R}_3(T)$ and ${\cal
R}_5(T)$ has height 
$$ 2\ell(T)=\frac{2}{\pi} (e^T + T) \sin\left(\frac{\pi T}{e^T +
T}\right), $$ 
and the larger horizontal strip of the regions ${\cal R}_1(T)$ and
${\cal R}_4(T)$ has height $2R(T)=\frac{2}{\pi} (e^T + T)$.
The arc of $\partial D_T$ cut out by the smaller strip has
length $2T$.

We construct an approximate solution on $Y(r(T))\times \R$ supported
on the strip $|t| \leq R(T)$.

We consider smooth functions $\eta_{T,1}(s,t)$, $\eta_{T,2}(s,t)$,
$\eta_{T,3}(s,t)$, $\eta_{T,4}(s,t)$, and $\eta_{T,5}(s,t)$ with
values in $[0,1]$, supported in the corresponding shaded regions in
Figure \ref{figII4}, respectively, such that  
$$\eta_{T,1} + \eta_{T,2} + \eta_{T,3}+\eta_{T,4}+\eta_{T,5} \equiv 1$$ 
is satisfied everywhere in (\ref{rectangle}).
We extend by $1$ the functions 
$$ \eta_{T,1}(s,t), \ \ \eta_{T,3}(s,t), \ \ \eta_{T,4}(s,t), \ \
\hbox{and} \ \ \eta_{T,5}(s,t)$$  
on the sides $V$ and $\nu(K)$ of the corresponding domains ${\cal
R}_i(T)$ which are not represented in the figure.  
 
We finally choose a smooth cutoff function $\chi(t)$ with values in
$[0,1]$, supported in $[-R(T),R(T)]$ and
satisfying $\chi(t)\equiv 1$ for $t\in [-R(T)+1, R(T)-1]$. Consider
also cutoff functions $\chi_\pm (t)$ supported in $[R(T)-1,\infty)$
and $(-\infty, -R(T)+1]$, and with $\chi_\pm (t)\equiv 1$ for $|t|\geq
R(T)$. 
We define the approximate solution as
\ba \begin{array}{ll} (\A,\Psi)_T=& \chi \cdot
(\eta_{T,1}(A(t),\psi(t))_T  + 
\eta_{T,2} a_T \\[2mm] & + \eta_{T,3} (\A',\Psi') + \eta_{T,4} (a ''(t),
0)_T +\eta_{T,5}(a_\infty '',0) )\\[2mm] & 
+ \chi_- (A^-,\psi^-)_T + \chi_+ (A^+,\psi^+)_T.\end{array}
\label{approx:sol} \na 
Here we use the notation
$(A(t),\psi(t))_T=(A(R(T)\tau),\psi(R(T)\tau))$, where
$(A(\tau),\psi(\tau))$ are the paths in $\M_V^*$, that appear in the
geometric limits, 
for $\tau\in [\ell(T)/R(T), 1]$, and $t=R(T)\tau$. We also denote by
$(A^\pm,\psi^\pm)_T$ representatives of the elements in $\M_{Y(r(T))}$
given by the gluing
$$ (A^\pm,\psi^\pm)_T = (A^\pm, \psi^\pm)\#_{r(T)} (a^\pm,0), $$
with
$$ (A^\pm, \psi^\pm)=(A(\pm 1),\psi(\pm 1)) \in \M_V^*. $$

\subsection{Linearizations}

Let $L_{A(t),\psi(t)}$ be the linearization of the 3-dim
Seiberg-Witten equations (\ref{3SW}) at the solutions
$(A(t),\psi(t))$. Consider 
the operator $L_{A(t),\psi(t)}$ acting on pairs $(\alpha,\phi)$ of
$L^2_1$ 1-forms and spinors.

First, let us recall from \cite{CMW} that, given an element
$[A,\psi]$ in ${\cal M}_V^*$, the operator $L_{A,\psi}$ is a Fredholm
operator on the weighted Sobolev spaces of 1-forms and spinors on $V$, 
\ba L_{A,\psi}: L^2_{1,\delta}(V)\to L^2_\delta(V), \label{Ldelta} \na
with 
$$\delta =\frac{1}{2} \min \{ | \lambda |, \ \ \lambda \in
spec(Q_{a_\infty'}) \}, $$
and $a_\infty'= \partial_\infty([A,\psi])$.
By Theorem 6.2 and 7.4 of \cite{LM}, we can consider the Fredholm
operator $L_{A,\psi}$ acting on
$L^2_1$ 1-forms and spinors. For a generic choice of the
perturbation of the equations on $V$, as in \cite{CMW}, the operator 
\ba L_{A,\psi}: L^2_1(V)\to L^2(V) \label{LL21} \na
is injective and surjective. In the following section we shall discuss
the effect of adding the perturbation. For the purpose of this section,
let us just assume that injectivity and surjectivity are achieved.
On the manifold $V_r$, with $r\geq r_0$, we can then consider
the operator $L_{A,\psi}$ acting between
\ba \label{LVT}
L_{A,\psi}: L^2_1(V_r, P_+ \oplus \ell)\to L^2(V_r), \na
that is, on the space of 1-forms and spinors on $V_r$, with boundary 
conditions defined by the APS condition \cite{APS}
$$ P_+ =L^2 \hbox{-closure of} \ \{ \phi_k | \mu_k >0 \}, $$
with $\phi_k$ eigenfunctions of the operator $Q_{a_\infty}$ on $T^2$
with eigenvalue $\mu_k$, and by the Lagrangian subspace $\ell$ in
$H^1(T^2,i\R)$ defined by the asymptotic values of the extended
$L^2$-solutions of $L_{A,\psi}(\alpha,\phi)=0$ on $V$.
By Proposition 2.4 of \cite{CLM}, we have that (\ref{LVT}) is a
self-adjoint Fredholm operator with
$$ Ker(L_{A,\psi}: L^2_1(V)\to L^2(V))\cong Ker(L_{A,\psi}: L^2_1(V_r,
P_+ \oplus \ell)\to L^2(V_r)). $$

\begin{Lem}
Consider the operators $L_{A(t),\psi(t)}$, $L_{A,\psi}$, and $L_{\tilde
A, \tilde \psi}$ acting on $L^2_1$ 1-forms and spinors on the
manifold $V$ with the infinite cylindrical end $T^2\times [0,\infty)$.
Assume that these operators $L_{A(t),\psi(t)}$ are surjective for all
$t\in [-1,0) \cup (0,1]$, and that the operators $L_{A,\psi}$ and
$L_{\tilde A, \tilde \psi}$ are also surjective. Moreover, assume also
that the path $[A(t),\psi(t)]$ and the elements $[A,\psi]$ and
$[\tilde A, \tilde \psi]$ in ${\cal M}_V$ are contained in the
complement of $\partial_\infty^{-1}(U_\vartheta)$, where
$\partial_\infty$ is the asymptotic value map, and $U_\vartheta$ is a
small neighborhood of the bad point $\vartheta$ in the character variety
of $T^2$. Then, for all $T >T_0$, the operator
$$ {\cal T}_T= R(T)^{-1} \frac{\partial}{\partial t} + L_{A(t),\psi(t)}  $$
acting on pairs $(\alpha,\phi)$ in the  $L^2_1$ completion of
compactly supported 1-forms and spinors in the domain 
\ba \label{S1} \begin{array}{c} {\cal S}_1(T)= (V\times
[-1,1])\cup_{T^2\times \{ s=0  
\}\times [-1,1]} 
\\[2mm] T^2\times \{ (s,t) \ | \  s\in [0,\tilde\Upsilon(t,T)], \ \
t\in [-1,1] \} \end{array} \na 
is also surjective. 
\label{Tsurj}
\end{Lem}

\proof 
From the previous discussion we can conclude that, assuming the
operators $L_{A(t),\psi(t)}$, $L_{A,\psi}$ and
$L_{\tilde A, \tilde \psi}$, with domain as in (\ref{LL21}) are surjective,
we obtain that the operator $L_{A(t),\psi(t)}$, with domain as in
(\ref{LVT}), with $r=\tilde\Upsilon(t,T)$, has trivial kernel and
is surjective. Similarly, the operators $L_{A,\psi}$ and
$L_{\tilde A, \tilde \psi}$, with domain as in (\ref{LVT}), with
$r=r_0$, have trivial kernel and are surjective. 

For simplicity of notation, we just write $L$ for the surjective
operators we are considering. Given any $h$, we
want to find a function $f$ in the domain satisfying  
$$ h = (\frac{1}{R(T)}\frac{\partial}{\partial t} + L )f. $$
We know that we have $h=L g$ for some $g$. Moreover, we have
an estimate 
$$ \| g(t) \|_{L^2(V_{\tilde\Upsilon(t,T)} \times \{ t \})} \leq C_t
\| h(t) \|_{L^2(V_{\tilde\Upsilon(t,T)} \times \{ t \})}, $$
which follows from Parseval's formula (cf. \cite{APS} pg.50),
where the constant $C_t$ is determined by the smallest absolute value of
the non-zero eigenvalues of the asymptotic operator $Q_{a_\infty(t)}$,
$$ C_t = \frac{c}{\lambda_\infty(t)}, $$
with 
$$ \lambda_\infty(t)=\min \{ | \lambda | \ | \ \lambda\in
spec(Q_{a_\infty(t)}) \}. $$

In the character variety $\chi(T^2)$ consider a small open
neighborhood $U_\vartheta$ of the ``bad point'' $\vartheta$ (see
\cite{CMW}). We have the asymptotic value map
$$\partial_\infty: {\cal M}_V \to \chi_0(T^2,V)$$
defined as in \cite{CMW}, with
$$ \partial_\infty(A(t),\psi(t))=a_\infty(t). $$
Notice that, if we assume that the path $(A(t),\psi(t))$, for $0< t
\leq 1$ is in the complement of $\partial_\infty^{-1}(U_\vartheta)$, then
the constant $C_t$ can be replaced with the uniform constant
$$ C =  c \cdot (\min_{0\leq t \leq 1} \lambda_\infty(t))^{-1}, $$
independent of $t$ and of the gluing parameter $T$.

So we can write
$$ h = (\frac{1}{R(T)}\frac{\partial}{\partial t} + L )g
-\frac{1}{R(T)} \frac{\partial}{\partial t} g. $$
Again by the surjectivity of $L$ we can write 
$$ -\frac{\partial}{\partial t} g = L g_1. $$
Again, we have an estimate $\| g_1 \|\leq C \| g \|$, hence iterating the
process we get
$$  h = \left(\frac{1}{R(T)}\frac{\partial}{\partial t} + L
\right)\left(\sum_{k=0}^n 
(-1)^k \left(\frac{1}{R(T)}\right)^k 
g_k\right) +(-1)^{n+1} \left(\frac{1}{R(T)}\right)^{n+1}
\frac{\partial}{\partial t} g_n. $$ 
For all $T\geq T_0$ we obtain 
$$ \sum_{k=0}^\infty  \left(\frac{1}{R(T)}\right)^k  C^k < \infty, $$
since $R(T)\sim e^T$ for large $T\geq T_0$.
Thus, as $n\to\infty$, we have uniform convergence of the series
$$  \sum_{k=0}^n (-1)^k \left(\frac{1}{R(T)}\right)^k g_k \to f, $$
and of the sequence
$$ (-1)^{n+1} \left(\frac{1}{R(T)}\right)^{n+1}
\frac{\partial}{\partial t} g_n \to 0. $$

\endproof

Over the unit disk $D$, consider the operator
$$ D_1=\left(\begin{array}{cc} \partial_\rho -* \partial_\theta &
R(T)^{-2} *d \\[2mm]
-e^{2\rho}*d & \partial_\rho \end{array}\right) $$
acting on pairs  $(\tilde a,\tilde f)$ where $\tilde a(w,\rho,\theta)$
is  1-form  and $\tilde f(w,\rho,\theta)$ is a function.
Consider also the operator
$$ D_2=\left(\begin{array}{cc} \frac{1}{R(T)}(\partial_\rho -i
\partial_\theta) & 
-e^{\rho+i\theta} i\bar\partial_a^* \\[3mm]
e^{\rho-i\theta} i\bar\partial_a &  \frac{1}{R(T)}(\partial_\rho +i
\partial_\theta) 
\end{array} \right), $$
acting on a spinor $(\tilde \alpha,\tilde \beta)$. Here
$a(\rho,\theta)$ is the image  
under the map $a: D\to \chi_0(T^2,Y)$ of Theorem \ref{limits},
restricted to a subdomain of $D$ homeomorphic to the disk. That is,
for each $(\rho,\theta)$, $a(\rho,\theta)$ is a flat connection on $T^2$.

\begin{Lem}
Without loss of generality, we can assume that the image of the map
$a: D\to \chi_0(T^2,Y)$ avoids the lattice of bad points $\vartheta$ in
$\chi_0(T^2,Y)$
\label{avoid}
\end{Lem}

\proof
We have unique bad point $\vartheta$ in the character variety of $T^2$,
\cite{CMW}, which corresponds to a lattice of bad points in the
universal cover $\chi_0(T^2,Y)$. 

First notice that, according to the proof of Lemma \ref{ts:limits},
the condition that the endpoint $a_\infty'$ is the bad point
$\vartheta$ is non-generic: by direct inspection of the characteristic
polynomial of the system (\ref{tildeM}) this case corresponds to the
vanishing of the coefficients $C_{nlk}$ that give the faster decay
$e^{\pm n\rho -C_{nlk} e^{\rho}}$. Thus, under generic choice of the
perturbation, as will be discussed in the next section, we can assume
that $a_\infty ' \neq \vartheta$. 

This implies that, under the maps $a_V$ and
$a_\nu$ of Theorem \ref{limits}, the image of the side
$\theta \in \{ -\pi/2, \pi/2 \}$ of the boundary of the half disk
$D^+$ is a smooth path in 
$\partial_\infty {\cal M}_V$ or in $\partial_\infty {\cal
M}_{\nu(K)}$, connecting the points $a^\pm$ which avoids the lattice of 
bad points. We can always replace the domain $D^+$ with
some smaller domain containing the same component $\theta \in \{
-\pi/2, \pi/2 \}$ and $\rho \in (-\infty, 0]$
of the boundary, in such a way that the image of the rest of the
domain also avoids
the lattice of bad points. That still determines uniquely the
holomorphic 
function. Thus, we can assume for simplicity that the image of the
subdomain of $D$ under the map $a: D\to \chi_0(T^2,Y)$
also contains no bad point $\vartheta$ in $\chi_0(T^2,Y)$.

\endproof

We want to consider $D_1 \oplus D_2$ acting on elements $(\tilde
a,\tilde f,\tilde \alpha,\tilde \beta)$ that represent deformations of 
$(a,0,0,0)$ in the $L^2_1$ completion of compactly supported 1-forms
and spinors on $D$.

\begin{Lem}
The operator $D_1\oplus D_2$ acting on the $L^2_1$ completion of
compactly supported 1-forms and spinors on $D$, has trivial Cokernel.
\label{Tsurj2}
\end{Lem}

\proof
We know \cite{CMW} that the operator 
$$ \left(\begin{array}{cc}0  & -e^{\rho+i\theta} \bar\partial_a^*
\\[2mm]  
e^{\rho-i\theta} \bar\partial_a & 0 \end{array} \right) $$
is surjective if the flat connection $a(\rho,\theta)$ avoids the
lattice of bad points $\vartheta$ in $\chi_0(T^2,Y)$. Then, by the same
argument of Lemma \ref{Tsurj}, we obtain that for sufficiently large
$T\geq T_0$, the operator $D_2$ is surjective.

Suppose given a pair $(\beta,h)$ of a 1-form and a
function, such that $(\beta,h)$  is orthogonal to the range of $D_1$.
Then, the pair $(\beta,h)$ satisfies the equations
$$ \partial_\rho h = R(T)^{-2} d^*\beta $$
$$ \partial_\rho \beta = *\partial_\theta \beta + *d (e^{-2\rho}h). $$

The Fourier coefficients satisfy the ODE
\ba \begin{array}{l}
h_{nlk}'=(R(T))^{-2}(-il v_{nlk} + in u_{nlk} ) \\[2mm]
u_{nlk}'=in v_{nlk} -e^{-2\rho} ik h_{nlk} \\[2mm]
v_{nlk}'=-in u_{nlk} + e^{-2\rho} il h_{nlk}. \end{array}
\label{coker:ODE} \na

We can show that solutions of (\ref{coker:ODE}) diverge sufficiently
fast at $\rho\to -\infty$, so that they fail to be in $L^2$. 
Since we are interested in the $\rho\to -\infty$ behavior, we can
isolate in the system (\ref{coker:ODE}) a leading term given by the
system 
\ba \begin{array}{l}
h_{nlk}'=R(T)^{-2}(-il v_{nlk} + in u_{nlk} ) \\[2mm] 
u_{nlk}'=-e^{-2\rho} ik h_{nlk}\\[2mm]
v_{nlk}'= e^{-2\rho} il h_{nlk}, \end{array} \label{coker:leading} \na
and treat the remaining terms as a perturbation.

Solutions of the system (\ref{coker:leading}) are solutions of 
$$ h_{nlk}''= e^{-2\rho}(R(T))^{-2} (k^2+l^2) h_{nlk}. $$
These are combinations of Bessel functions of type $I(0,z)$ and
$K(0,z)$, with the variable $z=c_{nlk} \sqrt{(R(T))^{-2} e^{-2\rho}}$,
for some constants $c_{nlk}>0$. None of these functions is in
$L^2(D)$.

\endproof

Consider the operator ${\cal T}_{a''(t)}$ of the form
\ba \frac{1}{R(T)} \frac{\partial}{\partial t} +
\left(\begin{array}{ccc} *d & 
-d & 0 \\ -d^* & 0 & 0 \\ 0 & 0 & \partial_{a''(t)} \end{array}
\right),  \label{D4} \na
where $\partial_{a''(t)}$ is the 3-dimensional Dirac operator, twisted
with the flat connection  
$$ a''(w,s,t)=a''(t)(w,s). $$ 

\begin{Lem}
If the path $a''(t)$ of connections in $\chi_0(T^2,Y)$ avoids the
lattice of bad points, then the operator ${\cal T}_{a''(t)}$, acting on
the $L^2_1$ completion of the space of compactly supported 1-forms and
spinors on the domain 
\ba \label{S4} \begin{array}{c} {\cal S}_4(T)=(\nu(K)\times
[-1,1])\cup_{T^2\times \{ 
s=0 \}\times [-1,1]} \\[2mm] T^2\times \{ (s,t) \ | \  s\in
[0,\tilde\Upsilon(t,T)], \ \ t\in [-1,1] \} \end{array} \na 
is surjective for sufficiently large $T\geq T_0$.
\label{Tsurj3}
\end{Lem}

\proof
The operator $$\frac{1}{R(T)} \frac{\partial}{\partial t}
+\partial_{a''(t)}$$ is 
surjective, for sufficiently large $T\geq T_0$, by the argument of
Lemma \ref{Tsurj}. Suppose given an element $(\beta,h)$ orthogonal to
the range of  
$$ \frac{1}{R(T)} \frac{\partial}{\partial t} +\left(\begin{array}{cc}*d & 
-d\\ -d^* & 0  \end{array} \right). $$
In the large $T >> 0$ limit, the element $(\beta,h)$ satisfies
$$ d^*\beta =0 \ \ \hbox{and} \ \ *d\beta -dh=0. $$
This implies that $h$ is harmonic, in the $L^2_1$ completion of the
space of compactly supported function, hence it is vanishing by 
the maximum principle. On the other hand $\beta$ satisfies
$(d+d^*)\beta =0$, which also implies $\beta\equiv 0$ in the $L^2_1$
completion of the space of compactly supported forms on (\ref{S4}).

\endproof

\subsection{The solutions}

In this section we prove that every approximate solution can be
deformed to an actual solution of the Seiberg--Witten equations. We
first prove that the linearization ${\cal D}_T$
at the approximate solution is
surjective. Here the operator 
$$ {\cal D}_T= {\cal D}_{{\A,\Psi}}$$ 
is the linearization of the 4-dimensional Seiberg--Witten equations
(\ref{4SW}) on $Y(r(T))\times \R$, and $(\A,\Psi)$ is the approximate
solution (\ref{approx:sol}). Then, the gluing theorem follows as a
fixed point argument, , as we are going to
discuss, cf. similar arguments in \cite{CMW},
\cite{DoSa}, \cite{Fl}, \cite{MW}, \cite{Sch}. 
Our case here is similar to the construction of \cite{DoSa}, since we
want to prove a gluing theorem where the underlying geometry also
changes along with the gluing parameter $T$. 
In particular, we will obtain the actual solution close to the
approximate solution by a fixed point arguments with constants
depending on $T$. 

In particular, we have seen from the convergence result of Theorem
\ref{limits}, that the non-uniform geometric limits are obtained from
solutions on $Y(r)\times \R$ after suitable rescaling of the
coordinates. Thus, in order to prove the gluing theorem, we will
introduce a norm on the domain $Y(r(T))\times \R$ which takes into
account the necessary rescaling. That will allow us to compare the
operator ${\cal D}_T$ with the operators $\D_i^T$ defined in
the previous subsection
on the corresponding regions ${\cal S}_i(T)$,
$$  \begin{array}{c} {\cal S}_1(T)= (V\times
[-1,1])\cup_{T^2\times \{ s=0  
\}\times [-1,1]} 
\\[2mm] T^2\times \{ (s,t) \ | \  s\in [0,\tilde\Upsilon(t,T)], \ \
t\in [-1,1] \} \end{array} $$
${\cal S}_2(T) = D$, the unit disk. ${\cal S}_i(T)={\cal R}_i(T)$ for
$i=3,5$, and 
$$ \begin{array}{c} {\cal S}_4(T)=(\nu(K)\times
[-1,1])\cup_{T^2\times \{ 
s=0 \}\times [-1,1]} \\[2mm] T^2\times \{ (s,t) \ | \  s\in
[0,\tilde\Upsilon(t,T)], \ \ t\in [-1,1] \}. \end{array} $$

We introduce a $T$-dependent Banach norm, which takes into account
the effect of rescaling. We need to define rescaled $L^2$ and $L^2_1$
norms on the spaces $\Lambda^1 \oplus \Gamma(W^+)$ and $\Lambda^{2+}
\oplus \Gamma(W^-)$ on $\cup_i {\cal R}_i(T)$.
We introduce the notation
$$ \epsilon(T)=\frac{1}{R(T)}. $$

First notice the following rescaling of the volume elements and norms.
Let us denote by
$$ \varphi: {\cal S}_i(T) \to {\cal R}_i(T) $$
the diffeomorphism given by the rescaling. In particular, we analyze
the two distinct cases
$$ \varphi: {\cal S}_1(T) \to {\cal R}_1(T) $$
$$ \varphi(\tilde x, \tilde y, \tilde z, \tilde t)=(x=\tilde x,
y=\tilde y, z=\tilde z , t=R(T) \tilde t), $$
and
$$ \varphi: {\cal S}_2(T) \to {\cal R}_2(T) $$
$$ \varphi(\tilde u, \tilde v, \tilde s, \tilde t)=(u=\tilde u,
v=\tilde v, s=R(T) \tilde s, t=R(T) \tilde t). $$
We have
$$ \int_{{\cal R}_i(T)} \omega = \int_{{\cal S}_i(T)} \varphi^*(\omega),
$$
for any 4-form $\omega$. Moreover, the Hodge $*$-operator satisfies
$$ \tilde *_{{\cal S}_1(T)} e^i = e^j \wedge e^k \wedge e^t $$
$$ \tilde *_{{\cal S}_1(T)} e^t = e^i\wedge e^j \wedge e^k $$
where the orthonormal basis of 1-form in the rescaled metric on ${\cal
S}_1(T)$ is given by $e^i = dx^i$ and $e^t = R(T)
d\tilde t$.
Similarly, the Hodge $*$-operator satisfies the same relations for the
orthonormal basis of 1-forms $du$, $dv$, $R(T) d\tilde s$,
$R(T) d\tilde t$ on ${\cal S}_2(T)$.

Consider first the space $\Lambda^1\oplus \Gamma(W^+)$. We write a
1-form as $a_V + f dt$ and a spinor $\psi$ on ${\cal R}_1(T)$. We
denote by $\tilde a_V$, $\tilde f$, $\tilde \psi$, etc. for the
composite $f\circ \varphi$ with the change of coordinates.
We obtain
$$ \int_{{\cal R}_1(T)} a_V \wedge *_4 a_V =\int_{{\cal S}_1(T)}
\varphi^*(a_V \wedge *_4 a_V) =  \int_{{\cal S}_1(T)}
\varphi^*(a_V) \wedge \tilde *_4 \varphi^*(a_V), $$
$$ \int_{{\cal R}_1(T)} f\, dt \wedge *_4 (f\, dt) = \int_{{\cal
S}_1(T)} \varphi^*(f\, dt \wedge *_4 (f\, dt)) = $$
$$ \int_{{\cal S}_1(T)} (R(T) \tilde f d\tilde t) \wedge 
\tilde *_4  (R(T) \tilde f d\tilde t)= \int_{{\cal
S}_1(T)} \varphi^*(f\, dt) \wedge \tilde *_4 \varphi^*(f\, dt). $$
$$ \int_{{\cal R}_1(T)} \| \psi \|^2 dt= R(T) \int_{{\cal
S}_1(T)} \| \tilde \psi \|^2 d\tilde t. $$

Thus, over this region, the elements $\xi_T= (a_v, R(T) f\, dt,
\epsilon(T)^{1/2} \psi)$ satisfy 
$$ \| \xi_T \|_{L^2({\cal R}_1(T))} = \| \varphi^*(\xi) \|_{L^2({\cal
S}_1(T))}, $$
for  $\xi=(a_v,  f\, dt, \psi)$.

Now consider the case of $\Lambda^{2+}\oplus \Gamma(W^-)$, again over
the region ${\cal R}_1(T)$. We have 2-forms $\Omega= \omega_V + \eta_V \wedge
dt$ and a spinors $\psi$ satisfying
$$ \int_{{\cal R}_1(T)} \omega_V \wedge *_4 \omega_V= \int_{{\cal S}_1(T)}
\varphi^*(\omega_V \wedge *_4 \omega_V) = \int_{{\cal S}_1(T)}
\varphi^*(\omega_V) \wedge \tilde *_4 \varphi^*(\omega_V), $$
$$ \int_{{\cal R}_1(T)} (\eta_V \wedge
dt)\wedge *_4 (\eta_V \wedge dt) = \int_{{\cal S}_1(T)}\varphi^*(
(\eta_V \wedge dt)\wedge *_4 (\eta_V \wedge dt)) = $$
$$  \int_{{\cal S}_1(T)} \varphi^* (\eta_V \wedge dt)
\wedge \tilde *_4 \varphi^* (\eta_V \wedge dt), $$
$$ \int_{{\cal R}_1(T)} \| \psi \|^2 dt= R(T) \int_{{\cal
S}_1(T)} \| \tilde \psi \|^2 d\tilde t. $$

Thus, in this case again, the elements $\xi_T= ( \omega_V,
 \eta_V\wedge dt, 
\epsilon(T)^{1/2} \psi)$ satisfy 
$$ \| \xi_T \|_{L^2({\cal R}_1(T))} = \| \varphi^*(\xi) \|_{L^2({\cal
S}_1(T))}, $$
for  $\xi=(\omega_V,\eta_V\wedge dt,\psi)$.

Now consider the regions ${\cal R}_2(T)$ and ${\cal S}_2(T)$. In the
case of $\Lambda^1\oplus \Gamma(W^+)$, we have 1-forms $a_{T^2} + f \,
ds + h \, dt$ and spinors $\psi=(\alpha,\beta)$, satisfying
$$ \int_{{\cal R}_2(T)} a_{T^2} \wedge *_4 a_{T^2}= \int_{{\cal
S}_2(T)} \varphi^*(a_{T^2} \wedge 
*_4 a_{T^2}) = \int_{{\cal S}_2(T)}
\varphi^*(a_{T^2})\wedge \tilde *_4 \varphi^*(a_{T^2}), $$
$$ \int_{{\cal R}_2(T)} f\, ds\wedge *_4 (f\, ds) = \int_{{\cal
S}_2(T)} \varphi^*(f\, ds\wedge *_4 (f\, ds)) = $$
$$ \int_{{\cal S}_2(T)} \varphi^*(f\, ds) \wedge \tilde *_4
\varphi^*(f\, ds), $$
$$ \int_{{\cal R}_2(T)} h\, dt\wedge *_4 (h\, dt) = \int_{{\cal
S}_2(T)} \varphi^*(h\, dt\wedge *_4 (h\, dt)) = $$
$$ \int_{{\cal S}_2(T)} \varphi^*(h\, dt) \wedge \tilde *_4
\varphi^*(h\, dt), $$
$$ \int_{{\cal R}_2(T)} \| \psi \|^2 ds\,dt= R(T)^2 \int_{{\cal
S}_2(T)} \| \tilde \psi \|^2 d\tilde s\, d\tilde t. $$

Thus, the elements $\xi_T=( a_{T^2}, f\, ds,h\, dt,
\epsilon(T)\psi)$ satisfy 
$$ \| \xi_T \|_{L^2({\cal R}_2(T))} = \| \varphi^*(\xi) \|_{L^2({\cal
S}_2(T))}$$
for $\xi=(a_{T^2}, f\, ds,h\, dt,  \psi)$   

Similarly, for $\Lambda^{2+}\oplus \Gamma(W^-)$ we have 2-forms
$\omega_{T^2}+ \eta\wedge ds + \gamma \wedge dt + h \, ds\wedge dt$,
and spinors $\psi$ satisfying
$$ \int_{{\cal R}_2(T)} \omega_{T^2} \wedge *_4
\omega_{T^2}=\int_{{\cal
S}_2(T)} \varphi^*(\omega_{T^2} \wedge *_4
\omega_{T^2}) =  \int_{{\cal S}_2(T)} \varphi^*(\omega_{T^2})
\wedge \tilde *_4 \varphi^*(\omega_{T^2}), $$
$$ \int_{{\cal R}_2(T)} \eta\wedge ds \wedge *_4 (\eta\wedge ds) =
 \int_{{\cal S}_2(T)}\varphi^*( \eta\wedge ds)\wedge \tilde *_4
\varphi^*( \eta\wedge ds), $$
$$ \int_{{\cal R}_2(T)} \gamma \wedge dt \wedge *_4 (\gamma \wedge dt) =
 \int_{{\cal S}_2(T)}\varphi^*(\gamma \wedge dt )\wedge
\tilde *_4 \varphi^*(\gamma \wedge dt ), $$   
$$ \int_{{\cal R}_2(T)} h\, ds\wedge dt \wedge *_4 (h\, ds\wedge dt)
= \int_{{\cal S}_2(T)}  \varphi^*(h\, ds\wedge dt) \wedge
\tilde *_4 \varphi^*(h\, ds\wedge dt), $$
$$ \int_{{\cal R}_2(T)} \| \psi \|^2 ds\,dt= R(T)^2 \int_{{\cal
S}_2(T)} \| \tilde \psi \|^2 d\tilde s\, d\tilde t. $$

Thus, the elements $\xi_T=( \omega_{T^2}, \eta\wedge ds,
\gamma \wedge dt,  h, \epsilon(T)\psi)$ satisfies
$$ \| \xi_T \|_{L^2({\cal R}_2(T))} = \| \varphi^*(\xi) \|_{L^2({\cal
S}_2(T))}$$
for $\xi=(\omega_{T^2}, \eta\wedge ds,\gamma \wedge dt, h,\psi)$.

Now consider the following $T$-rescaled $L^2$-norm for an element
$\xi$ in $\Lambda^1\oplus \Gamma(W^+)$ over $\cup_i {\cal R}_i(T)$:
\ba \label{rescaled:norm:0} \begin{array}{ll} \| \xi \|^2_{0,T}:= &
\| (\eta_{T,1}+\eta_{T,4}) (a_V,f\, dt, \epsilon(T)^{1/2} \psi)
\|^2_{L^2({\cal R}_1(T) \cup {\cal R}_4(T))} + \\[2mm]
& \| (\eta_{T,3}+\eta_{T,5}) (a_V,f\, dt, \psi) \|^2_{L^2({\cal R}_3(T)
\cup {\cal R}_5(T))} + \\[2mm] 
& \| \eta_{T,2} (a_{T^2}, f\, ds, h\, dt,
\epsilon(T)(\alpha,\beta))\|^2_{L^2({\cal R}_2(T))}. \end{array} \na
Similarly, we define the $T$-rescaled $L^2$-norm for an element
$\xi$ in $\Lambda^{2+}\oplus \Gamma(W^-)$ over $\cup_i {\cal R}_i(T)$: 
\ba \label{rescaled:norm:0:2-} \begin{array}{ll} \| \xi \|^2_{0,T}:= &
\| (\eta_{T,1}+\eta_{T,4}) (\omega_V,\eta_V\wedge dt, \epsilon(T)^{1/2} \psi)
\|^2_{L^2({\cal R}_1(T) \cup {\cal R}_4(T))} + \\[2mm]
& \| (\eta_{T,3}+\eta_{T,5}) (\omega_V,\eta_V\wedge dt, \psi)
\|^2_{L^2({\cal R}_3(T) 
\cup {\cal R}_5(T))} + \\[2mm] 
& \| \eta_{T,2} (\omega_{T^2}, \eta\wedge ds, \gamma\wedge dt, h,
\epsilon(T)(\alpha,\beta))\|^2_{L^2({\cal R}_2(T))}. \end{array} \na
With this choice of norms we have 
$$ \| \xi \|_{0,T\, {\cal R}_i(T)} = \| \varphi^*(\xi) \|_{L^2({\cal
S}_i(T))}. $$

Moreover, we define a $T$-rescaled $L^2_1$ norm on $\Lambda^1\oplus
\Gamma(W^+)$ on $\cup_i {\cal R}_i(T)$ by setting
\ba \label{rescaled:norm:1} \begin{array}{ll} 
\| \xi \|^2_{1,T} : = & \| (\eta_{T,1}+\eta_{T,4}) (a_V,f\, dt,
\epsilon(T)^{1/2} \psi)  \|^2_{L^2 ({\cal R}_1(T) \cup {\cal R}_4(T))}
+ \\[2mm] & \sum_{i=1}^3 \| (\eta_{T,1}+\eta_{T,4}) \nabla_i (a_V,f\, dt,
\epsilon(T)^{1/2} \psi)  \|^2_{L^2 ({\cal R}_1(T) \cup {\cal R}_4(T))}
+ \\[2mm] & \| (\eta_{T,1}+\eta_{T,4}) R(T)\frac{\d}{\d t} (a_V,f\, dt,
\epsilon(T)^{1/2} \psi)  \|^2_{L^2 ({\cal R}_1(T) \cup {\cal R}_4(T))}
+ \\[2mm] & 
\| (\eta_{T,3}+\eta_{T,5}) (a_V,f\, dt, \psi) \|^2_{L^2({\cal R}_3(T)
\cup {\cal R}_5(T))} + \\[2mm] 
& \sum_{i=1}^3 \| (\eta_{T,3}+\eta_{T,5}) (a_V,f\, dt, \psi)
\|^2_{L^2({\cal R}_3(T) 
\cup {\cal R}_5(T))} + \\[2mm] 
& \| (\eta_{T,3}+\eta_{T,5}) R(T) \frac{\d}{\d t} (a_V,f\, dt, \psi)
\|^2_{L^2({\cal R}_3(T)\cup {\cal R}_5(T))} + \\[2mm] 
& \| \eta_{T,2} (\omega_{T^2}, \eta\wedge ds, \gamma\wedge dt, h,
\epsilon(T)(\alpha,\beta))\|^2_{L^2({\cal R}_2(T))} + \\[2mm] 
& \sum_{i=1}^3 \| \eta_{T,2} \nabla_i (\omega_{T^2}, \eta\wedge ds,
\gamma\wedge dt, h, 
\epsilon(T)(\alpha,\beta))\|^2_{L^2({\cal R}_2(T))} + \\[2mm] 
&  \| \eta_{T,2} R(T) \frac{\d}{\d t} (\omega_{T^2}, \eta\wedge ds,
\gamma\wedge dt, h, 
\epsilon(T)(\alpha,\beta))\|^2_{L^2({\cal R}_2(T))} \end{array} \na

On $\Lambda^1\oplus
\Gamma(W^+)$ on the ${\cal S}_i(T)$ we consider the 
$L^2_1$ norm defined as
$$ \begin{array}{ll} \| \xi \|^2_{1,T\, {\cal S}_1(T)\cup {\cal
S}_3(T)}: = & 
\| \xi \|^2_{L^2({\cal S}_1(T)\cup {\cal S}_3(T))} + \\[2mm]  &
\sum_{i=1}^3 \| \nabla_i \xi \|^2_{L^2({\cal S}_1(T)\cup {\cal
S}_3(T))} + \|  \frac{\d}{\d \tilde t}  \xi \|^2_{L^2({\cal
S}_1(T)\cup {\cal S}_3(T))}, \end{array} $$
$$ \begin{array}{ll} \| \xi \|^2_{1,T\, {\cal S}_3(T)\cup {\cal
S}_5(T)}: = & 
\| \xi \|^2_{L^2({\cal S}_1(T)\cup {\cal S}_3(T))} + \\[2mm]  &
\sum_{i=1}^3 \| \nabla_i \xi \|^2_{L^2({\cal S}_1(T)\cup {\cal
S}_3(T))} + \|  \frac{\d}{\d t}  \xi \|^2_{L^2({\cal
S}_1(T)\cup {\cal S}_3(T))}, \end{array} $$
$$ \begin{array}{ll} \| \xi \|^2_{1,T\, {\cal S}_2(T)} : = & \| \xi
\|^2_{L^2({\cal S}_2(T))} + \\[2mm]  &
\| \d_u \xi \|^2_{L^2({\cal S}_2(T))} +  \| \d_v \xi
\|^2_{L^2({\cal S}_2(T))} + \\[2mm]  &
\| \d_{\tilde t} \xi \|^2_{L^2({\cal S}_2(T))} +  \| 
\d_{\tilde s} \xi \|^2_{L^2({\cal S}_2(T))}. \end{array} $$
With this choice of norms, we have $$ \| \xi \|_{1,T\, {\cal R}_i(T)}= \|
\varphi^*(\xi) \|_{L^2_1 {\cal S}_i(T)}. $$

We can now prove the main Lemma for the gluing theorem: this can be
regarded as an analogue, in our context, of Lemma 4.4 and 4.5 of
\cite{DoSa}. We denote in the following by $\tilde \D_i^{T}$ the
operators $\D_i^{T}$ introduced before, rescaled under
${\cal S}_i(T) \to {\cal R}_i(T)$. 

\begin{Lem} Suppose given a sequence $\xi_k$ of connections and
spinors, which we can write as 
$$ \xi_k = (A_k, f_k, \Psi_k) $$
over the domains $\cup_{i\neq 2}{\cal R}_i(T_k)$, and 
$$ \xi_k = (a_k,f_k,h_k,\alpha_k,\beta_k) $$
over ${\cal R}_2(T_k)$, for a sequence $T_k\to
\infty$.
Let $\xi_{k,i}=\eta_{T_k,i} \xi_k$. 
We assume that the elements $\xi_{k,i}$ are in the orthogonal
complement of the Kernels of
the operators $\tilde \D_i^{T_k}$, 
$$ \xi_{k,i} \in Ker(\tilde \D_i^{T_k})^\perp, $$
with respect to the $L^2_{1,T_k}$-norms.
Moreover, we assume that the operators $\tilde \D_i^{T_k}$ have
trivial cokernels in the space $L^2_{0,T_k}({\cal R}_i(T_k))$. 
Then, under this hypothesis, the convergence
$$ \| {\cal D}_{T_k} \xi_k \|_{0,T_k} \to 0, $$
where ${\cal D}_{T_k}$ is the linearization at the approximate
solutions, implies
$$ \| \xi_k \|_{1,T_k} \to 0. $$
\label{glue:1}
\end{Lem}

\proof We first observe that, in the operator norm over each ${\cal
R}_i(T_k)$, we have 
$$ \| \D_{T_k,i} - \tilde\D_i^{T_k} \| \to 0 $$
as $T_k\to \infty$. 
We also assume that the cutoff functions satisfy $$sup_{{\cal
R}_i(T_k)} | \nabla \eta_{T_k,i}(s,t) | \leq q(T_k).$$ We assume
that $q(T_k) \to 0$, where a bound on the rate of decay to zero of
$q(T)$ as $T\to \infty$ will be specified in the proof of Lemma
\ref{glue:lem:1} below. 

The argument then is similar to \cite{Sch} \S 2.5. 
Suppose given a sequence $\xi_k$,  
and parameters $T_k\to\infty$, satisfying 
$$ \| {\cal D}_{T_k} \xi_k \|_{0,T_k} \to 0. $$ Moreover, assume that,
for all $k$, we have
$$ \| \xi_k \|_{1,T_k} = 1. $$
We have an estimate
$$ \| \tilde {\cal D}_i^{T_k} \xi_{k,i} \|_{0,T_k} \leq C q(T_k) \| \xi_k
\|_{1,T_k}+ \|  \eta_{T_k,i} \tilde {\cal D}_i^{T_k} \xi_k  \|_{0,T_k}
$$ 
$$ \leq C q(T_k) + \| {\cal D}_{T_k} - \tilde {\cal D}_i^{T_k} \| \,
\| \xi_k \|_{1,T_k} + \| {\cal D}_{T_k} \xi_k \|_{0,T_k}. $$ 
The terms in the right hand side decay to zero, thus we obtain that
the elements $\xi_{k,i}$ are in the span of the low modes of $\tilde
{\cal D}_i^{T_k}$. Let us denote
by ${\cal V}_{k,i}$ the space of eigenvectors of $(\tilde {\cal
D}_i^{T_k})^* \tilde {\cal D}_i^{T_k}$ acting on completion of the
space of compactly supported 1-forms and spinors over ${\cal
R}_i(T_k)$ in the $L^2_{1,T_k}$ norm, with eigenvalues $\lambda_{T_k}
\to 0$ as $T_k\to \infty$. This is the space of low modes of $\tilde
{\cal D}_i^{T_k}$. 
We have obtained $\xi_{k,i} \in {\cal V}_{k,i}$, from the previous
estimate. We use the notation ${\cal V}_{k,i}^{\#}$ for the span of
low modes of the $L^2_{1,T_k}$ adjoint.  

Now we have 
$$ \dim {\cal V}_{k,i} \geq \dim Ker(\tilde {\cal D}_i^{T_k}). $$
Moreover, for $T_k\to \infty$, we have
$$ \dim {\cal V}_{k,i} - \dim {\cal V}_{k,i}^{\#} = \dim Ker(\tilde
{\cal D}_i^{T_k}) - \dim Coker(\tilde {\cal D}_i^{T_k}). $$
Under the assumption that $Coker(\tilde {\cal D}_i^{T_k})=0 $, we
obtain the reverse estimate
$$ \dim {\cal V}_{k,i} \leq \dim Ker(\tilde {\cal D}_i^{T_k}). $$
Thus, we can identify ${\cal V}_{k,i} \cong Ker(\tilde {\cal
D}_i^{T_k})$. Thus, under these hypotheses we would have $\xi_{k,i} \in
Ker(\tilde {\cal D}_i^{T_k})$. 
This contradicts the initial assumption $\xi_{k,i} \in Ker(\tilde
{\cal D}_i^{T_k})^\perp$. 

\endproof

The assumption $Coker(\tilde {\cal D}_i^{T_k})=0 $ has been discussed
in the previous subsection, and it follows from the results of the
last Section, about perturbed operators.

We can derive from Lemma \ref{glue:1} the following Corollary.

\begin{Cor}
There are constants $T_0$ and $c_1>0$, independent of $T$, such that we
have an estimate 
$$ \| \xi \|_{1,T} \leq c_1  \| {\cal D}_T \xi \|_{0,T}, $$ 
for all $\xi$ satisfying 
$$ \eta_{T,i} \xi \in Ker(\tilde {\cal D}_i^T)^\perp. $$
Moreover, we
have a similar estimate
$$ \| \xi \|_{0,T} \leq C  \| {\cal D}^*_T \xi \|_{1,T}, $$
for a constant $C >0$ independent of $T\geq T_0$, 
under the assumption that the operators $\tilde \D_i^T$ have trivial
Cokernels. Here ${\cal D}^*_T$ is the adjoint with respect to the
$T$-dependent norm. \label{glue:1:c} 
\end{Cor}

\proof By the choice of the $T$-dependent norms,  $c_1$ is independent
of $T$. The rest of the statement
follows from Lemma \ref{glue:1}. The second estimate follows by a
similar argument.  

\endproof

Now we can state the second main result of this part of the work,
namely the gluing theorem.

\begin{The}
Given any approximate solution $\Xi_0=(\A,\Psi)$, constructed as in
(\ref{approx:sol}), for all sufficiently large $T\geq T_0$, there
exists a solution $\Xi$ of the equations (\ref{4SW}) on $Y(r(T))\times
\R$, satisfying $\| \Xi - \Xi_0 \|_{1,T}\leq c \epsilon(T)^{1/2}$,
for a constant $c>0$ independent of $T\geq T_0$.
\label{glue:theorem}
\end{The}

\proof We divide the proof in several steps.
On the domain $Y(r(T))\times \R$, consider the map
$$ \sigma_T (\A,\Psi)=\left\{ \begin{array}{l}
F^+_{\A}-\tau(\Psi,\Psi) \\[2mm] D_{\A}\Psi \end{array}\right. $$ 

\begin{Lem} Let $\Xi_0$ be the approximate solution as in
(\ref{approx:sol}). The estimate
$$ \| \sigma_T (\Xi_0) \|_{0,T} \leq c_0 \epsilon(T)^{1/2} $$
is satisfied, with $c_0 >0$ independent of $T\geq T_0$. 
\label{glue:lem:1}
\end{Lem}

\proof
Consider first the region ${\cal R}_2(T)\cup{\cal R}_3(T)\cup {\cal
R}_5(T)$. Over this region we have 
$\sigma_T (\Xi_0)(w,s,t) =0$ except on the overlap between the
supports of the functions $\eta_{T,i}$ used in
(\ref{approx:sol}). 
On ${\cal R}_i(T)\cap supp(\eta_{T,j})$ in this region, we can
estimate
$$ \| \sigma_T (\Xi_0) \|_{0,T\, ({\cal R}_i(T)\cap supp(\eta_{T,j}))}
\leq \sup_{{\cal R}_j(T)} |\nabla \eta_{T,j}| \, \, \| (\A,\Psi) \|_{0,T \,
({\cal R}_i(T))}. $$
Recall that we have assumed $\sup_{{\cal R}_j(T)} |\nabla \eta_{T,j}|=
q(T)$. It is sufficient to
choose the cutoff functions as in Lemma \ref{glue:1}, with
the hypothesis that $q(T) \sim \epsilon(T)^{1/2}$, and we get the desired
estimate.  Then consider the case of the region ${\cal 
R}_1(T)$. Here the condition $\sigma_T (\Xi_0)(w,s,t) =0$ is not
satisfied, but we can estimate the error term by
$$ \| \sigma_T (\Xi_0) \|^2_{0,T\, ({\cal R}_1(T))} \leq 
\int_{I(T)} ( \| \frac{\partial}{\partial t} A(t)
\|_{L^2(V)}^2 + \epsilon(T) \| \frac{\partial}{\partial t} \psi(t)
\|^2_{L^2(V)}) dt $$
$$ \leq \epsilon(T) \int_{-1}^1 \| \frac{\partial}{\partial \tilde t}
(A(\tilde t),\psi(\tilde t)) \|^2_{L^2(V)} d\tilde t \leq c^2
\epsilon(T). $$ 
The remaining case of the region
${\cal R}_4(T)$ is analogous, with a similar resulting estimate with
the constant
$$ c \geq \left(\int_{-1}^1 \| \partial_{\tilde t} a''(\tilde t)
\|^2_{L^2(T^2)} d\tilde t \right)^{1/2}, $$
where $a''(\tilde t)$ is the path in ${\cal M}_{\nu(K)}$ obtained in the
geometric limits.

\endproof

Now define elements as follows
$$ \Xi_1= \Xi_0 + \xi_0 \ \ \ \xi_0 =\D_T^* \eta_0 \ \ \ \D_T \D_T^*
\eta_0= -\sigma_T(\Xi_0). $$
The elements are well defined because of the following.

\begin{Lem}
The equation $\D_T \D_T^*\eta_0= -\sigma_T(\Xi_0)$ admits a unique
solution $\eta_0$ with $\| \eta \|_{0,T}\leq C$. Moreover, we have
$$ \| \xi_0 \|_{1,T} \leq c_1 \| \sigma_T(\Xi_0)
\|_{0,T}\leq c_0 c_1 \epsilon(T)^{1/2}.  $$ 
with $c_1 >0$ independent of $T\geq T_0$. 
\label{glue:lem:2}
\end{Lem}

\proof The operator $\D_T \D_T^*$ is invertible, under the assumptions
of Lemma \ref{glue:1}. The estimate then follows
from the previous Lemma \ref{glue:1}, Corollary \ref{glue:1:c} and Lemma 
\ref{glue:lem:1}. 

\endproof

The linearization of the map $\sigma_T$ at $\Xi_0=(\A,\Psi)$ is given by
$$ d\sigma_{T,\Xi_0}(\eta, \Phi) = \left\{ \begin{array}{l}  
d^+ \eta -\frac{1}{2} Im(\Psi, \Phi) \\[2mm]
D_{\A} \Phi + \eta . \Psi. \end{array}\right. $$  
By construction we have
$$ d\sigma_{T,\Xi_0}(\xi)= \D_T(\xi). $$

Consider the non-linear part of the map $\sigma_T$, that is the expression
$$ {\cal N}\sigma_{T}(\xi)=\sigma_T(\Xi_0
+\xi)-\sigma_T(\Xi_0)-d\sigma_{T,\Xi_0}(\xi). $$

We have
$$\sigma_T(\Xi_1)=\sigma_T(\Xi_0+\xi_0) -\sigma_T(\Xi_0)
-d\sigma_{T,\Xi_0}(\xi_0)={\cal N}\sigma_{T}(\xi_0). $$

\begin{Lem}
We have an estimate
$$ \| \sigma_T(\Xi_1) \|_{0,T} \leq c_2  \| \xi_0 \|_{0,T}^2
\leq c_0 c_1 c_2 \epsilon(T)^{1/2} \| \xi_0 \|_{0,T}, $$
\label{glue:lem:3}
\end{Lem}

\proof
The non-linear part is given by 
$$ {\cal N}\sigma_{T}(\Omega, \Phi) =\left\{ \begin{array}{l}
\tau(\Phi,\Phi) \\[2mm] \Omega . \Phi, \end{array}\right. $$ 
for $(\Omega, \Phi)$ in $\Lambda^1\oplus \Gamma(W^+)$. Thus, in the
$L^2$-norms we have 
$$ \| {\cal N}\sigma_{T}(\xi_0) \|_{L^2} \leq c_2
 \| \xi_0 \|^2_{L^2}. $$
This is sufficient to obtain the desired estimate on the regions
${\cal R}_3(T)\cup {\cal R}_5(T)$. 
If we write $\Omega= \omega + f dt$, we can write the above estimate
more precisely as 
$$ \| \omega . \Phi \|_{L^2} + \| f dt . \Phi \|_{L^2} \leq c (\|
\omega  \|_{L^2} \| \Phi \|_{L^2} + \|  f \|_{L^2} \| \Phi
\|_{L^2}), $$
and similarly, for the term $\tau(\Phi,\Phi)$, we have  $\|
\tau(\Phi,\Phi) \|_{L^2} \leq c \| \Phi \|_{L^2}^2$.
We have
$$ \| \tau(\Phi,\Phi) \|_{0,T}= \|
\varphi^*(\tau(\Phi,\Phi)) \|_{L^2({\cal S}_i(T))}, $$
and
$$ \| \Omega .\Phi \|_{0,T}=  \|
 \varphi^*(\Omega .\Phi)\|_{L^2({\cal S}_i(T))}, $$
by our choice of weights on $\Lambda^{2+}\oplus \Gamma(W^-)$.
Moreover, we have $$\varphi^*(\tau(\Phi,\Phi))=\tilde\tau(\tilde \Phi,
\tilde \Phi) \ \ \hbox{ and } \ \ \varphi^*(\Omega
.\Phi)=\varphi^*(\Omega)\tilde . \tilde \Phi. $$ Thus, we obtain
$$ \| \tau(\Phi,\Phi) \|_{0,T} \leq c_2  \| \tilde \Phi
\|^2_{L^2({\cal S}_i(T))} $$
$$ \leq c_2  \| \Phi \|^2_{0,T \,({\cal
R}_i(T))},  $$ 
where the last inequality comes form the rescaling of the norm on
$\Lambda^1\oplus \Gamma(W^+)$.
Similarly, we have 
$$ \| \Omega .\Phi\|_{0,T} \leq c_2 
\|\varphi^*(\Omega) \|_{L^2({\cal S}_i(T))}
 \cdot \| \tilde \Phi \|_{L^2({\cal S}_i(T))} $$
$$ \leq c_2  \|\Omega \|_{0,T \,({\cal R}_i(T))}\cdot \|
\Phi \|_{0,T \,({\cal R}_i(T))}. $$
Thus, we have obtained an estimate
$$ \| \sigma_T(\Xi_1) \|_{0,T}\leq  c_2  \| \xi_0 \|_{0,T}^2. $$
We also have an estimate
$$ \| \xi_0 \|_{1,T} \leq c_0 c_1 \epsilon(T)^{1/2}, $$
from the previous Lemmas, which gives 
$$ \| \sigma_T(\Xi_1) \|_{0,T}\leq  c_0 c_1 c_2 \epsilon(T)^{1/2} \|
\xi_0 \|_{0,T}. $$ 
This proves the claim.

\endproof

We then proceed inductively, as in \cite{DoSa}.
We set
$$ \Xi_{\nu +1}=\Xi_\nu + \xi_\nu \ \ \ \xi_\nu=\D_T^* \eta_\nu \ \ \
\D_T \D_T^* \eta_\nu= -\sigma_T(\Xi_\nu). $$

\begin{Lem}
The following estimates hold:
$$ \| \xi_\nu \|_{1,T} \leq C(\nu) \epsilon(T)^{1/2} $$ 
$$ \|  \sigma_T(\Xi_{\nu+1}) \|_{0,T} \leq \hat C(\nu)  \| \xi_\nu
\|_{0,T}. $$ 
Moreover, we can always assume that an estimate
$$ \| \sigma_T(\Xi_1) \|_{0,T}\leq  c_0 c_1 c_2 \|
\xi_0 \|_{0,T}, $$ 
with $c=c_0 c_1^2 c_2 < 1$, is satisfied. This implies that the
iteration process 
converges to a solution $\Xi$ of the SW equations on $Y(r(T))\times
\R$, satisfying $\| \Xi-
\Xi_0 \|_{1,T} \leq C \epsilon(T)^{1/2}$.
\label{glue:lem:4}
\end{Lem}

\proof First of all notice that we can always include a factor
$\epsilon(T_0)^{1/2}$ in the constant $c_0 c_1 c_2$ in the estimate of
Lemma \ref{glue:lem:3}, for sufficiently large $T_0$, so that the
condition $c_0 c_1^2 c_2 < 1$ is satisfied.
Then define recursively $C(0)=c_0 c_1$ and $\hat C(0)= c_0 c_1 c_2$ and
$C(\nu+1)=c_0 c_1^2 c_2 C(\nu)$ and $\hat C(\nu +1) = c_0 c_1^2 c_2
\hat C(\nu)$. The result then follows inductively using the 
estimates
$$ \| \xi_\nu \|_{1,T} \leq c_1  \| \sigma_T(\Xi_\nu)
\|_{0,T}, $$
as in Lemma \ref{glue:1} and
$$ \| \sigma_T(\Xi_{\nu +1 }) \|_{0,T}\leq c_2  \| \xi_\nu
\|_{1,T}^2, $$ as in Lemma \ref{glue:lem:3}.
The estimate for the solution $\Xi$ is obtained from
$$ \| \Xi_\nu -\Xi_0 \|_{1,T} \leq \sum_{j=0}^{\nu -1} \| \xi_j
\|_{1,T}. $$

\endproof

This completes the proof of the gluing Theorem \ref{glue:theorem}.

\endproof

\section{Perturbation}

The gluing theorem stated in the previous section 
relies on the assumption that 
the linearization $L_{A(t),\psi(t)}$ is surjective, for all the elements
$[A(t),\psi(t)]$ in ${\cal M}_V$, as in Lemma
\ref{Tsurj}. We know that, in general, these conditions are satisfied
only after introducing a suitable perturbation. 

In \cite{CMW} we defined a class of perturbations ${\cal P}$ of the
Chern--Simons--Dirac functional, such that, for a generic element
$P=(U,V)\in {\cal P}$, the operators $L_{A,\psi}$ and ${\cal
D}_{\A,\Psi}$ are both surjective. These operators, in the perturbed
case, are the linearizations of the
corresponding perturbed  critical point equation
$$ \left\{ \begin{array}{l}
*F_A = \sigma (\psi, \psi) +\displaystyle{
 \sum_{j=1}^{N} \frac{\d U}{\d \tau_j}}\mu_j \\[2mm]
\dirac_A (\psi)+ \displaystyle{\sum_{j=1}^{K} \frac{\d V}{\d \zeta_j}}
\nu_j.\psi =0.
\end{array}\right. $$
and of the perturbed flow line equation
$$ \left\{ \begin{array}{lll}
\displaystyle{\frac{\d A}{\d t} }&=& -*F_A + \sigma (\psi, \psi) + 
\sum_{j=1}^{N} \displaystyle{\frac{\d U}{\d \tau_j}} \mu_j \\[2mm]
\displaystyle{\frac{\d \psi}{\d t}}& =& -\dirac_A \psi 
- \sum_{j=1}^{K} \displaystyle{\frac{\d V}{\d \zeta_j}}\nu_j.\psi,
\end{array}\right. $$
respectively. 

For simplicity, we use the notation
$$ P_U = \sum_{j=1}^{N} \frac{\d U}{\d \tau_j}\mu_j, $$
$$ P_V = \sum_{j=1}^{K} \frac{\d V}{\d \zeta_j}\nu_j. $$
Thus, we have
\ba 
\left\{ \begin{array}{l}
*F_A = \sigma (\psi, \psi) + P_U \\[2mm]
\dirac_A (\psi)+ P_V . \psi 
\end{array}\right. \label{3SWP} \na
and
\ba \label{4SWP} \left\{ \begin{array}{l}
\frac{\d A}{\d t}=-*F_A + \sigma (\psi, \psi) + P_U \\[2mm]
\frac{\d \psi}{\d t}= -\dirac_A \psi -P_V . \psi.
\end{array}\right. \na

On the manifold $Y(r)$ with a long cylinder, and on the non-compact
manifold $V$ with an infinite cylinder, we refine the definition of
the class of perturbations, as in \cite{CMW}, by including the
requirement that the 
perturbation is exponentially small in the region  $T^2\times [-r,r]
\cup \nu(K)$, and on the end $T^2\times [0,\infty)$. We denote the
corresponding class with ${\cal P}_\delta$, where $\delta$ is the rate 
of decay, depending on the smallest absolute value of the non-trivial
eigenvalues of the asymptotic operator $Q_{a_\infty}$ on $T^2$,
cf. \cite{CMW}. 
 
We need to check that the main results in the previous sections can be 
extended to the case where the perturbation $P\in {\cal P}_\delta$ is
introduced. 

\subsection{Estimates in the perturbed case}

In order to extend the results of Section 2 to the case of the
equations (\ref{4SWP}), we need the analogue of the uniform estimate
on the energy, and then of Lemma \ref{L1} and
Lemma \ref{Lpoint}. The results of Section 2 then extend to this case
without significant changes.

\begin{Lem}
Let $(\A_r,\Psi_r)$ be a finite energy solution of (\ref{4SWP}) on
$Y(r)\times \R$. Then, for $r\geq r_0$, and for any interval
$[t_0,t_1]$ of length $\ell=t_1-t_0$, the estimates of Lemma \ref{L1}
and Lemma \ref{Lpoint} hold, with
$$ s_0 =\max_{Y(r_0)} \{ -s(x) +C(P), 0 \}, $$
where $C(P)$ is a positive constant depending only on the perturbation 
$P \in {\cal P}_\delta$.
\label{L1P}
\end{Lem}

\proof Recall that the class of perturbations ${\cal P}_\delta$ is
defined \cite{CMW} by considering complete $L^2$ bases $\{ \nu_i(r)
\}_{i=1}^\infty$ and $\{ \mu_j(r) \}_{j=1}^\infty$, satisfying
$$ \sup_{T^2\times [-r,r]} | \nu_i(r) | \leq \sup_{T^2\times
[-r_0,r_0]} | \nu_i(r_0) | $$ 
$$ \sup_{T^2\times
[-r,r]} | \mu_j(r) | \leq \sup_{T^2\times [-r_0,r_0]} | \mu_j(r_0) |, $$
and rescaling them with a function
$$ f_r(s)=e^{-\delta (s+r)}, $$
for $-r+\epsilon \leq s \leq r-\epsilon$, on the cylinder $T^2\times
[-r,r]$, with the weight $\delta$ satisfying
\[
\delta \geq \frac 12 min \{ \lambda_{a_\infty} | a_\infty \in
\chi(T^2)\backslash U_\vartheta \}.
\]
The elements $\{ f_r \nu_i(r) \}_{i=1}^\infty$ and $\{ f_r \mu_j(r)
\}_{j=1}^\infty$ still give complete bases, which we use to define the 
perturbation $P_r$ in the class ${\cal P}_\delta$ on $Y(r)$.

Thus, the pointwise estimate obtained from the
Weitzenb\"ock formula gives
\ba \begin{array}{c}
0 \geq \frac{s}{2}|\psi|^2 - \la *F_A. \psi, \psi\ra + \la \frac{\d
U}{\d\zeta_i} f_r\nu_i(r).\psi,  \frac{\d U}{\d\zeta_j}
f_r\nu_i(r).\psi\ra \\[2mm] 
 - \la \frac{\d U}{\d\zeta_i} d(f_r\nu_i(r)).\psi,\psi\ra 
\geq \frac{s}{2}|\psi|^2 +\frac 12 |\psi|^4 - C(P_r)
|\psi|^2. \end{array} \label{unifestP} \na
The constant satisfies $C(P_r)\leq C(P_0)$. 

We can prove a uniform
bound on the energy. In fact, the energy is now defined as the
variation of the perturbed Chern--Simons--Dirac functional
$$ CSD_P(A,\psi)= CSD (A, \psi) + U(\tau_1 (A, \psi), \cdots, \tau_N
(A, \psi)) $$
$$ + V (\zeta_1(A, \psi), \cdots, \zeta_K (A, \psi)), $$
as in \cite{CMW}. 

If $(\A_r,\Psi_r)$ is a family of finite energy solutions on
$Y(r)\times \R$, with asymptotic values
$(A_r(\pm\infty),\psi_r(\pm\infty))$ as $t\to\pm\infty$ satisfying the
perturbed equations (\ref{3SWP}), the energy
$$ \E_r = CSD_P(
A_r(-\infty),\psi_r(-\infty))-CSD_P(A_r(+\infty),\psi_r(+\infty)) $$
is given by
$$ \frac{1}{2} \int_{Y(r)} F_{\A_r}\wedge F_{\A_r} + \int_{Y(r)} \la
\psi_r(-\infty), f_r \sum \frac{\d V_r}{\d \zeta_j}
\nu_j(r). \psi_r(-\infty) \ra $$
$$ -\int_{Y(r)} \la \psi_r(+\infty), f_r \sum \frac{\d V_r}{\d \zeta_j}
\nu_j(r). \psi_r(+\infty) \ra $$
$$ + U(\tau_1 (A_r(-\infty), \psi_r(-\infty)), \cdots, \tau_N
(A_r(-\infty), \psi_r(-\infty))) $$
$$ - U(\tau_1 (A_r(+\infty),
\psi_r(+\infty)), \cdots, \tau_N (A_r(+\infty), \psi_r(+\infty))) $$
$$ + V(\zeta_1 (A_r(-\infty), \psi_r(-\infty)), \cdots, \zeta_K
(A_r(-\infty), \psi_r(-\infty))) $$
$$ - V(\zeta_1 (A_r(+\infty),
\psi_r(+\infty)), \cdots, \zeta_K (A_r(+\infty), \psi_r(+\infty))). $$
We know the first term is uniformly bounded, from the analysis of the
unperturbed case. Using the notation 
$$ P_V(r)=f_r \sum \frac{\d V_r}{\d \zeta_j}
\nu_j(r), $$ 
the second and third term can be estimated by
writing
\ba \begin{array}{c}
| \int_{Y(r)} \la \psi_r, P_V(r). \psi_r \ra dv | \leq \int_{Y(r_0)} | \la
\psi, P_V(r). \psi \ra | dv \\[2mm] + \int_{T^2\times ([-r,-r_0]\cup [r_0,r])}
| \la \psi_r, P_V(r). \psi_r \ra | dv. \end{array} \label{1est} \na
The first term in the right hand side is bounded by
$$ \| \sigma(\psi_r, \psi_r) \|_{L^2(Y_(r_0))} \cdot \| P_V(r)
\|_{L^2(Y(r_0))} \leq C\ C(P_0)^2 \| P_V(r_0) \|_{L^2(Y(r_0))}, $$
where the terms on the right are obtained using the uniform pointwise
estimate (\ref{unifestP}).
The second term in the right hand side of (\ref{1est}) can be
estimated similarly by
$$ C (r-r_0)^2 \  C(P_0)^2  e^{-\delta (r-r_0)} \sup | \nu_i(r_0) |, $$
where the factor $(r-r_0)^2$ comes from factoring out the volume of the
cylinder in the estimate of both $\sigma(\psi_r,\psi_r)$ and
$P_V(r)$. Again we have used the pointwise estimate (\ref{unifestP}). 

The remaining terms in the variation of $CSD_P$ can be bounded as
follows. We have
$$ | U(\tau_j(A_r(-\infty),\psi_r(-\infty))) -
U(\tau_j(A_r(+\infty),\psi_r(+\infty))) | $$
$$ \leq | \tau_j(A_r(-\infty),\psi_r(-\infty)) -
\tau_j(A_r(+\infty),\psi_r(+\infty)) | \cdot $$
$$ \| \frac{\d U(r)}{\d \tau_j} \mu_j(r) \|_{L^2(Y_r)}. $$
The last term is bounded uniformly, by our assumptions on the
perturbation. We have
$$ \tau_j(A_r(-\infty),\psi_r(-\infty)) -
\tau_j(A_r(+\infty),\psi_r(+\infty)) $$
$$ =\int_{Y(r)} f_r (A_r(-\infty) - A_r(+\infty))\wedge *\mu_j(r). $$
As in Lemma 4.6 of \cite{CMW}, up to changing the connections within
the same gauge class, we have an estimate
$$ \| A_r(-\infty) - A_r(+\infty) \|_{L^2(Y(r))} \leq C(P_0)\ Cr, $$
where the right hand side grows linearly in $r$ like the volume
$Vol(Y(r))$. 
This estimate follows from the uniform pointwise bound on the spinor, and the
corresponding bound on the curvature. Thus, for all $\epsilon >0$ we
can choose $r\geq r_0$ large enough so that we have a bound
$$ | \int_{Y(r)} f_r (A_r(-\infty) - A_r(+\infty))\wedge *\mu_j(r) | $$
$$ \leq C C(P_0) Vol(Y(r_0)) + \epsilon. $$
The estimate of the remaining term in the variation of $CSD_P$ is
analogous. Combining these estimates, we get a uniform bound on the
energy for large enough $r\geq r_0$, hence the estimates on the finite
energy solutions $(\A_r,\Psi_r)$ follow as in the Lemmata \ref{L1} and
\ref{Lpoint}. 

\endproof

\subsection{Asymptotics}

In order to adapt the results of Section 3  we need to study the
asymptotics of finite energy solutions of the equations
(\ref{4SWP}) on the manifold $V\times \R$, with the non-compact end
$T^2\times [0,\infty)\times \R$.

As in Section 3.2, we consider the ODE associated to the perturbed
system (\ref{4SWP}) on $V\times \R$.
We write (\ref{4SWP}) as 
$$ \partial_t a -dh + *(\partial_s a -df)= * i(\bar\alpha \beta
+\alpha\bar\beta)+ \sum_{j=1}^N \frac{\d U}{\d \tau_j} 
* q_j $$ 
$$ \partial_t f -\partial_s h + * F_a
=\frac{i}{2}(|\alpha|^2-|\beta|^2)+ \sum_{j=1}^N \frac{\d U}{\d
\tau_j}p_j $$ 
$$ \partial_t \alpha + h\alpha + i \partial_s \alpha + if \alpha +
\bar\partial_a^*\beta -i\sum_{j=1}^N \frac{\d V}{\d
\zeta_j}((-\nu_j^1+i\nu_j^2) \beta + i \nu_j^0 \alpha)=0 $$
$$ \partial_t \beta + h\beta -i \partial_s \beta -if\beta +
\bar\partial_a \alpha +i\sum_{j=1}^N \frac{\d V}{\d
\zeta_j}((\nu_j^1 +i\nu_j^2) \alpha - i \nu_j^0 \beta)=0, $$
where we use the notation introduced in \cite{CMW},
$$ *_3 P_U = \sum_{j=1}^N \frac{\d U}{\d \tau_j} (p_j + q_j\wedge ds) $$
and
$$ P_V=\sum_{i=1}^K \frac{\d V}{\d \zeta_i}(\nu_i^1 dx +\nu_i^2 dy
+\nu_i^0 ds), $$ 
with $*\mu_j=p_j + q_j\wedge ds$, and $\nu_i=\nu_i^1 dx +\nu_i^2 dy
+\nu_i^0 ds$. 
In radial gauge we obtain  
\ba
\label{SWradialP}
\begin{array}{l}
\partial_\rho f =e^{2\rho} *( F_a +
\frac{i}{2}(|\alpha|^2-|\beta|^2)\omega +P_0) \\[2mm]
\partial_\rho a =*(\partial_\theta a -df + i(\bar\alpha \beta
+\alpha\bar\beta)+P_1) \\[2mm]
\partial_\rho \alpha = i( \partial_\theta \alpha + f \alpha + e^{\rho
+i\theta} (\bar\partial_a^*\beta +P_{11}\alpha + P_{12}\beta) ) \\[2mm]
\partial_\rho \beta = -i( \partial_\theta \beta + f\beta +
e^{\rho-i\theta} (\bar\partial_a \alpha +P_{21}\alpha + P_{22}\beta),
\end{array} 
\na
where we use the notation
\ba \begin{array}{c} 
P_0(s)=\sum_{j=1}^N \frac{\d U}{\d \tau_j}p_j(s) \\[2mm]
P_1^{0,1}(s)=\sum_{j=1}^N \frac{\d U}{\d \tau_j} 
q_j^{0,1}(s) \\[2mm]
P(s)=\left(\begin{array}{cc} P_{11} & P_{12} \\ P_{21} & P_{22}
\end{array}\right) =
\sum_{j=1}^N \frac{\d V}{\d
\zeta_j} \left(\begin{array}{cc} i\nu_j^0 & -\nu_j^1 +i \nu_j^2 \\
\nu_j^1 +i \nu_j^2 & -i\nu_j^0 \end{array}\right) 
\end{array} \label{P}\na
as in \cite{CMW}

The perturbation terms $P_0(s)$, $P_1(s)$ and $P(s)$ depend on
$(a,f,\alpha,\beta)$ through the variables
$$ \tau_j= \int (A-A_0) \wedge * \mu_j $$
and
$$ \zeta_j=\int \la \nu_j. \psi, \psi \ra dv $$
(cf. \cite{CMW}). By the condition that the perturbation is chosen in
the class ${\cal P}_\delta$ we obtain the estimates
$$ \begin{array}{lr}
\| P_0(\rho) \|_{L^2(T^2\times \{\theta \})} \leq & C(U,V) \|
(a,f,h) \|_{L^2(T^2\times \{\theta \})} \\[2mm]
&  \cdot \int_{-\infty}^\infty \exp(-\delta e^\rho \cos \theta) d\rho
\end{array} $$
$$ \begin{array}{lr}
\| P_1(\rho) \|_{L^2(T^2\times \{\theta \})}\leq & C(U,V) \|
(a,f,h) \|_{L^2(T^2\times \{\theta \})} \\[2mm]
&  \cdot \int_{-\infty}^\infty \exp(-\delta e^\rho \cos \theta) d\rho 
\end{array} $$
\ba \begin{array}{lr}
\| P(\rho) \cdot (\alpha,\beta) \|_{L^2(T^2\times \{\theta \})}  \leq
& C(U,V) \|  
(\alpha,\beta)\|_{L^2(T^2\times \{\theta \})} \\[2mm]
& \cdot \int_{-\infty}^\infty \exp(-\delta e^\rho \cos \theta) d\rho
\end{array} \label{est:P} \na
after the change of coordinates $s+it = e^{\rho+i\theta}$.

In order to study the asymptotics of the system (\ref{SWradialP}) we
proceed as in Section 3. We consider the linear system given by the
uncoupled systems (\ref{ODE:F}) and (\ref{uncoupled}) for the linear
ASD and Dirac equations. We add the perturbation terms coming from
the terms (\ref{P}) in polar coordinates $(\rho,\theta)$.
We first study the asymptotic of this perturbed system and show that
the finite energy solutions are still exponentially decaying in the 
radial direction, as the solutions of the original systems
(\ref{ODE:F}) and (\ref{uncoupled}). We then proceed as in the
remaining of Section 3, to prove that the full system
(\ref{SWradialP}) has finite energy solutions that are exponentially
decaying in the radial direction. 

The system given by (\ref{ODE:F}) and (\ref{uncoupled}) together with
the perturbation terms (\ref{P}) is also uncoupled in the curvature
and Dirac part. We discuss the behavior of the Dirac part: the
curvature part is completely analogous.

Following \cite{Ha}, \S X, Section 8, we consider all the systems of
ODE's of the form (\ref{uncoupled}), for all $(n,l,k)$, with the
additional terms coming from the Fourier transform of the term
$P(\rho)\cdot (\alpha,\beta)$. These may no longer be uncoupled as the 
original (\ref{uncoupled}). We can write this perturbed system in the
form 
\ba \begin{array}{c}
X^\prime = M^- X + P^- (X,Y) \\
Y^\prime = M^+ Y + P^+ (X,Y), \end{array} \label{2system} \na
where the variables $X$ correspond to the eigenvectors of the systems
(\ref{uncoupled}) with eigenvalues $\lambda =\pm\lambda_i^{nlk}$, as
in (\ref{eigenv}) with $Re(\lambda) < 0$, and $Y$ corresponds to the
eigenvectors with $Re(\lambda) >0$.

Let $\lambda_0$ be the smallest absolute value of the eigenvalues with
$Re(\lambda) < 0$. Consider a fixed $\mu$ with $-\lambda_0< \mu < 0$.
According to Theorem 8.1 and 8.3 of \cite{Ha}, \S X, the finite energy 
solutions of (\ref{2system}) will be of the form
$$ X(\rho)= e^{M^-(\rho-\rho_0)} X_0 + \int_{\rho_0}^\rho
e^{M^-(\rho-\tau)}P^- (X(\tau),Y_0) d\tau $$
(cf. (8.15) of \cite{Ha}, \S X).

Thus, the asymptotic decay of solutions is governed by the decay of
solutions (\ref{U}) of (\ref{uncoupled}). The asymptotics of the
original system  (\ref{SWradialP}) are then obtained by successive
approximation as in Section 3. 

With these results in place, the remaining of Section 3 and Section 4 
extend with minor changes. The hypotheses of  Lemma \ref{Tsurj} are
now satisfied for a generic choice of 
the perturbation $P\in {\cal P}_\delta$. Similarly, we have guaranteed
that the linearization $\D_{\A',\Psi'}$ is surjective.

\vskip .3in

\noindent {\bf Matilde Marcolli}, Department of Mathematics,
Massachusetts Institute of Technology, 2-275. Cambridge, MA 02139, USA
\par
\noindent matilde\@@math.mit.edu

\vskip .2in

\noindent {\bf Bai-Ling Wang}, Department of Pure Mathematics,
University of Adelaide, Adelaide SA 5005, Australia. \par
\noindent bwang\@@maths.adelaide.edu.au

\end{document}